\documentclass{article}
\usepackage{amssymb}
\usepackage{amsfonts}
\usepackage{amsmath}

\newtheorem{theorem}{Theorem}

\newtheorem{corollary}[theorem]{Corollary}

\newtheorem{lemma}[theorem]{Lemma}

\newtheorem{proposition}[theorem]{Proposition}
\newtheorem{remark}[theorem]{Remark}

\input{tcilatex}
\begin{document}

\author{George Avalos\thanks{%
email address: gavalos@math.unl.edu. The research of G. Avalos was partially
supported by the NSF grants DMS-1211232 and DMS-1616425.} \\
Department of Mathematics\\
University of Nebraska-Lincoln, 68588 USA \and Pelin G. Geredeli\thanks{%
email address: pguvengeredeli2@unl.edu. The research of P.G. Geredeli was
partially supported by the NSF grant DMS-1616425 and an Edith T. Hitz
Fellowship.} \\
Department of Mathematics\\
University of Nebraska-Lincoln, 68588 USA\\
\\
and\\
Department of Mathematics\\
Hacettepe University, 06800 Beytepe, Ankara, Turkey }
\date{December 6, 2017}
\title{Exponential Stability and Supporting Spectral Analysis of a
Linearized Compressible Flow-Structure PDE Model}
\maketitle

\begin{abstract}
In this work, a result of exponential stability is obtained for solutions of
a compressible flow-structure partial differential equation (PDE) model
which has recently appeared in the literature. In particular, a compressible
flow PDE and its associated state equation for the associated pressure
variable, each evolving within a three dimensional domain $\mathcal{O}$, are
coupled to a fourth order plate equation which holds on a flat portion $%
\Omega $ of the boundary $\partial \mathcal{O}$. Moreover, since this
coupled PDE model is the result of a linearization of the compressible
Navier-Stokes equations about an arbitrary state, the flow PDE component
contains a generally nonzero ambient flow profile $\mathbf{U}$. By way of
obtaining the aforesaid exponential stability, a \textquotedblleft frequency
domain\textquotedblright\ approach is adopted here, an approach which is
predicated on obtaining a uniform estimate on the resolvent of the
associated flow-structure semigroup generator.
\end{abstract}

\section{Introduction}

\bigskip

\subsection{Description of the problem}

\subsection{The PDE Model}

In this work, we will consider spectral and uniform decay properties of a
coupled partial differential equation (PDE) system, a system in which the
distinct PDE dynamics each evolve on different domains and moreover contact
each other only through a boundary interface.

Throughout, the flow domain $\mathcal{O}$ will be a bounded subset of $%
\mathbb{R}^{3}$\thinspace , with boundary $\partial \mathcal{O}$. Moreover, $%
\partial \mathcal{O}=\overline{S}\cup \overline{\Omega }$, with $S\cap
\Omega =\emptyset $, and with (structure) domain $\Omega \subset \mathbb{R}%
^{3}$ being a \emph{flat} portion of $\partial \mathcal{O}$. In particular, $%
\partial \mathcal{O}$ has the following specific configuration: 
\begin{equation}
\Omega \subset \left\{ x=(x_{1,}x_{2},0)\right\} \,\text{\ and \ surface }%
S\subset \left\{ x=(x_{1,}x_{2},x_{3}):x_{3}\leq 0\right\} \,.  \label{geo}
\end{equation}%
So if $\nu (x)$ denotes the unit normal vector to $\partial \mathcal{O}$,
pointing outward, then 
\begin{equation}
\nu |_{\Omega }=\left[ 0,0,1\right] \,.  \label{normal}
\end{equation}%
In addition, $\left[ \mathcal{O},\Omega \right] $ is assumed to fall within
one of the following classes:%
\begin{equation*}
\begin{array}{l}
\mathsf{(G.1)}\ \mathcal{O}\text{ is a convex domain with wedge angles }<%
\frac{2\pi }{3}\text{. \ Moreover, }\Omega  \\ 
\text{ \ \ \ \ \ \ \ has Lipschitz boundary;} \\ 
\mathsf{(G.2)}\ \mathcal{O}\text{ is a convex polyhedron having angles}<%
\frac{2\pi }{3}\text{,} \\ 
\text{ \ \ \ \ \ \ \ \ \ \ \ \ \ \ and so then }\Omega \text{ is a convex
polygon with angles}<\frac{2\pi }{3}\text{.}%
\end{array}%
\end{equation*}

\medskip 

In addition, if $\ \mathbf{n}(\mathbf{x})$ denotes the unit outward normal
vector to $\partial \mathcal{O}$, then necessarily $\left. \mathbf{n}%
\right\vert _{\Omega }=(0,0,1)$

The coupled PDE system which we will consider is the result of a
linearization which is undertaken in \cite{Chu2013-comp} and \cite{agw}:
Within the three-dimensional geometry $\mathcal{O}$, the compressible
Navier-Stokes equations are present (see e.g., \cite{temam}), assuming the
flow which they describe to be barotropic.This system is linearized with
respect to some reference rest state of the form $\left\{ p_{\ast },\mathbf{U%
},\varrho _{\ast }\right\} $: the pressure and density components ${p_{\ast
},\varrho _{\ast }}$ are scalars, and the arbitrary ambient field $\mathbf{U}%
:\mathcal{O}\rightarrow \mathbb{R}^{3}$ is given here by: 
\begin{equation}
\mathbf{U}%
(x_{1},x_{2},x_{3})=[U_{1}(x_{1},x_{2},x_{3}),U_{2}(x_{1},x_{2},x_{3}),U_{3}(x_{1},x_{2},x_{3})].
\label{flowfield}
\end{equation}%
In \cite{Chu2013-comp} and \cite{agw} non-critical lower order terms are
deleted, and the aforesaid pressure and density reference constants are set
equal to unity. Thus, we are presented with the following system of
equations, in solution variables $\mathbf{u}(x_{1},x_{2},x_{3},t)$ (flow
velocity), $p(x_{1},x_{2},x_{3},t)$ (pressure), $w_{t}(x_{1},x_{2},t)$
(elastic plate displacement) and $w_{t}(x_{1},x_{2},t)$ (elastic plate
velocity): 
\begin{align}
& \left\{ 
\begin{array}{l}
p_{t}+\mathbf{U}\cdot \nabla p+div~\mathbf{u}=0~\text{ in }~\mathcal{O}%
\times (0,\infty ) \\ 
\mathbf{u}_{t}+\mathbf{U}\cdot \nabla \mathbf{u}-div~\sigma (\mathbf{u}%
)+\eta \mathbf{u}+\nabla p=0~\text{ in }~\mathcal{O}\times (0,\infty ) \\ 
(\sigma (\mathbf{u})\mathbf{n}-p\mathbf{n})\cdot \boldsymbol{\tau }=0~\text{
on }~\partial \mathcal{O}\times (0,\infty ) \\ 
\mathbf{u}\cdot \mathbf{n}=0~\text{ on }~S\times (0,\infty ) \\ 
\mathbf{u}\cdot \mathbf{n}=w_{t}~\text{ on }~\Omega \times (0,\infty )%
\end{array}%
\right.  \label{system1} \\
&  \notag \\
& \left\{ 
\begin{array}{l}
w_{tt}+\Delta ^{2}w+\left[ 2\nu \partial _{x_{3}}(\mathbf{u})_{3}+\lambda 
\text{div}(\mathbf{u})-p\right] _{\Omega }=0~\text{ on }~\Omega \times
(0,\infty ) \\ 
w=\frac{\partial w}{\partial \nu }=0~\text{ on }~\partial \Omega \times
(0,\infty )%
\end{array}%
\right.  \label{IM2} \\
&  \notag \\
& 
\begin{array}{c}
\left[ p(0),\mathbf{u}(0),w(0),w_{t}(0)\right] =\left[ p_{0},\mathbf{u}%
_{0},w_{0},w_{1}\right] .%
\end{array}
\label{IC_2}
\end{align}%
This flow-structure system is a generalization of that considered by the
late Igor Chueshov in \cite{Igor-note}; however, unlike the PDE system in 
\cite{Igor-note}, the PDE system (\ref{system1})-(\ref{IC_2}) depends upon a
generally \emph{non-zero}, fixed, ambient vector field $\mathbf{U}$ about
which the linearization takes place. Here (pointwise in time), $p(t):\mathbb{%
R}^{3}\rightarrow \mathbb{R}$ and $\mathbf{u}(t):\mathbb{R}^{3}\rightarrow 
\mathbb{R}^{3}$ are given as the pressure and the flow velocity field,
respectively. The quantity $\eta >0$ represents a drag force of the domain
on the viscous flow. In addition, the quantity $\mathbf{\tau }$ in (\ref%
{system1}) is in the space $TH^{1/2}(\partial \mathcal{O)}$ of tangential
vector fields of Sobolev index 1/2; that is,%
\begin{equation}
\mathbf{\tau }\in TH^{1/2}(\partial \mathcal{O)=}\{\mathbf{v}\in \mathbf{H}^{%
\frac{1}{2}}(\partial \mathcal{O})~:~\left. \mathbf{v}\right\vert _{\partial 
\mathcal{O}}\cdot \mathbf{n}=0~\text{ on }~\partial \mathcal{O}\}.\footnote{%
See e.g., p.846 of \cite{buffa2}.}  \label{TH}
\end{equation}

In addition, we define the space 
\begin{equation}
\mathbf{V}_{0}=\{\mathbf{v}\in \mathbf{H}^{1}(\mathcal{O})~:~\left. \mathbf{v%
}\right\vert _{\partial \mathcal{O}}\cdot \mathbf{n}=0~\text{ on }~\partial 
\mathcal{O}\}.  \label{V_0}
\end{equation}

\medskip

Moreover, the \textit{stress tensor} $\sigma $ in the flow PDE component of (%
\ref{system1}) is defined as, 
\begin{equation*}
\sigma (\mathbf{\mu })=2\nu \epsilon (\mathbf{\mu })+\lambda \lbrack
I_{3}\cdot \epsilon (\mathbf{\mu })]I_{3},
\end{equation*}%
where \textit{Lam\'{e} Coefficients }$\lambda \geq 0$ and $\nu >0$; and the 
\textit{strain tensor }$\epsilon $ is given by 
\begin{equation*}
\epsilon _{ij}(\mu )=\dfrac{1}{2}\left( \frac{\partial \mu _{j}}{\partial
x_{i}}+\frac{\partial \mu _{i}}{\partial x_{j}}\right) \text{, \ }1\leq
i,j\leq 3
\end{equation*}%
(see p. 129 of \cite[p.129]{kesavan}).

As we noted in \cite{agw}, the flow PDE boundary conditions which are in (%
\ref{system1}) are the so-called \emph{impermeability}-slip conditions \cite%
{bolotin,chorin-marsden}: namely, no flow passes through the boundary -- in
particular, the normal component of the flow field $\mathbf{u}$ on the
active boundary portion $\Omega $ matches the plate velocity $w_{t}$ -- and
on $\partial \mathcal{O}$ there is no stress in the tangential direction $%
\tau $.

\subsection{Notation}

Throughout, for a given domain $D$, the norm of corresponding space $%
L^{2}(D) $ will be denoted as $||\cdot ||_{D}$ (or simply $||\cdot ||$ when
the context is clear). Inner products in $L^{2}(\mathcal{O})$ or $\mathbf{L}%
^{2}(\mathcal{O})$ will be denoted by $(\cdot ,\cdot )_{\mathcal{O}}$,
whereas inner products $L^{2}(\partial \mathcal{O})$ will be written as $%
\langle \cdot ,\cdot \rangle _{\partial \mathcal{O}}$. We will also denote
pertinent duality pairings as $\left\langle \cdot ,\cdot \right\rangle
_{X\times X^{\prime }}$ for a given Hilbert space $X$. The space $H^{s}(D)$
will denote the Sobolev space of order $s$, defined on a domain $D$; $%
H_{0}^{s}(D) $ will denote the closure of $C_{0}^{\infty }(D)$ in the $%
H^{s}(D)$-norm $\Vert \cdot \Vert _{H^{s}(D)}$. We make use of the standard
notation for the boundary trace of functions defined on $\mathcal{O}$, which
are sufficently smooth: i.e., for a scalar function $\phi \in H^{s}(\mathcal{%
O})$, $\frac{1}{2}<s<\frac{3}{2}$, $\gamma (\phi )=\phi \big|_{\partial 
\mathcal{O}},$ which is a well-defined and surjective mapping on this range
of $s$, owing to the Sobolev Trace Theorem on Lipschitz domains (see e.g., 
\cite{necas}, or Theorem 3.38 of \cite{Mc}).

\section{Preliminaries}

\bigskip\ Our result of uniform decay for solutions of the compressible
flow-structure PDE model -- and the supporting spectral analysis on $i%
\mathbb{R}$ -- will be stated below within the context of the associated
semigroup formulation which was given in \cite{agw} to solve the coupled PDE
system (\ref{system1})--(\ref{IC_2}) with given finite energy initial data $%
\left[ p_{0},\mathbf{u}_{0},w_{0},w_{1}\right] $. In particular, the
associated space of well-posedness is 
\begin{equation}
\mathcal{H}\equiv L^{2}(\mathcal{O})\times \mathbf{L}^{2}(\mathcal{O})\times
H_{0}^{2}(\Omega )\times L^{2}(\Omega ).  \label{H}
\end{equation}%
$\mathcal{H}$ is a Hilbert space, topologized by the following inner
product: 
\begin{equation}
(\mathbf{y}_{1},\mathbf{y}_{2})_{\mathcal{H}}=(p_{1},p_{2})_{L^{2}(\mathcal{O%
})}+(\mathbf{u}_{1},\mathbf{u}_{2})_{\mathbf{L}^{2}(\mathcal{O})}+(\Delta
w_{1},\Delta w_{2})_{L^{2}(\Omega )}+(v_{1},v_{2})_{L^{2}(\Omega )}
\label{innerp}
\end{equation}%
for any $\mathbf{y}_{i}=(p_{i},\mathbf{u}_{i},w_{i},v_{i})\in \mathcal{H}%
,~i=1,2.$

In \cite{agw}, we established that solutions to the compressible
flow-structure PDE system (\ref{system1})--(\ref{IC_2}), with initial data
in said finite energy space $\mathcal{H}$, can be associated with a strongly
continuous semigroup $\left\{ e^{\mathcal{A}t}\right\} \,_{t\geq 0}\subset 
\mathcal{L}(\mathcal{H)}$, where the generator $\mathcal{A}:D(\mathcal{A}%
)\subset \mathcal{H}\rightarrow \mathcal{H}$ will be described below. That
is, if for given Cauchy data $\left[ p_{0},\mathbf{u}_{0},w_{0},w_{1}\right]
\in \mathcal{H}$, vector $\left[ p,\mathbf{u},w,w_{t}\right] \in C([0,T];%
\mathcal{H})$ solves the problem (\ref{system1})--(\ref{IC_2}), then it also
solves the following abstract ODE (and vice versa): 
\begin{eqnarray}
\dfrac{d}{dt}%
\begin{bmatrix}
p \\ 
\mathbf{u} \\ 
w \\ 
w_{t}%
\end{bmatrix}
&=&\mathcal{A}%
\begin{bmatrix}
p \\ 
\mathbf{u} \\ 
w \\ 
w_{t}%
\end{bmatrix}%
;  \notag \\
\lbrack p(0),\mathbf{u}(0),w(0),w_{t}(0)] &=&[p_{0},\mathbf{u}%
_{0},w_{0},w_{1}].  \label{ODE}
\end{eqnarray}

\noindent

\bigskip

In fact, in \cite{agw} it is shown that the abstract operator $\mathcal{A}:D(%
\mathcal{A})\subset \mathcal{H}\rightarrow \mathcal{H}$ which mathematically
models the flow-structure dynamics (\ref{system1})--(\ref{IC_2}) is given
explicitly by

\begin{equation}
\mathcal{A}=\left[ 
\begin{array}{cccc}
-\mathbf{U}\mathbb{\cdot }\nabla (\cdot ) & -\text{div}(\cdot ) & 0 & 0 \\ 
-\mathbb{\nabla (\cdot )} & \text{div}\sigma (\cdot )-\eta I-\mathbf{U}%
\mathbb{\cdot \nabla (\cdot )} & 0 & 0 \\ 
0 & 0 & 0 & I \\ 
\left. \left[ \cdot \right] \right\vert _{\Omega } & -\left[ 2\nu \partial
_{x_{3}}(\cdot )_{3}+\lambda \text{div}(\cdot )\right] _{\Omega } & -\Delta
^{2} & 0%
\end{array}%
\right] .  \label{AAA}
\end{equation}

\noindent In addition, it is shown in \cite{agw} that semigroup generation
ensues when the domain $D(\mathcal{A})$ is specified as follows:%
\begin{equation*}
D(\mathcal{A})=\{(p_{0},\mathbf{u}_{0},w_{1},w_{2})\in L^{2}(\mathcal{O}%
)\times \mathbf{H}^{1}(\mathcal{O})\times H_{0}^{2}(\Omega )\times
H_{0}^{2}(\Omega )~:~\text{properties }(A.i)\text{--}(A.v)~~\text{hold}\},
\end{equation*}%
where

\begin{enumerate}
\item[(A.i)] $\mathbf{U}\cdot \nabla p_{0}\in L^{2}(\mathcal{O}).$

\item[(A.ii)] $\text{div}~\sigma (\mathbf{u}_{0})-\nabla p_{0}\in \mathbf{L}%
^{2}(\mathcal{O})$. (Consequently, we infer the boundary trace regulatity 
\newline
$\left[ \sigma (\mathbf{u}_{0})\mathbf{n}-p_{0}\mathbf{n}\right] _{\partial 
\mathcal{O}}\in \mathbf{H}^{-\frac{1}{2}}(\partial \mathcal{O})$.)

\item[(A.iii)] $-\Delta ^{2}w_{0}-\left[ 2\nu \partial _{x_{3}}(\mathbf{u}%
_{0})_{3}+\lambda \text{div}(\mathbf{u}_{0})\right] _{\Omega }+\left.
p_{0}\right\vert _{\Omega }\in L^{2}(\Omega ).$

\item[(A.iv)] $\left( \sigma (\mathbf{u}_{0})\mathbf{n}-p_{0}\mathbf{n}%
\right) \bot ~TH^{1/2}(\partial \mathcal{O})$. That is, 
\begin{equation*}
\left\langle \sigma (\mathbf{u}_{0})\mathbf{n}-p_{0}\mathbf{n},\mathbf{\tau }%
\right\rangle _{\mathbf{H}^{-\frac{1}{2}}(\partial \mathcal{O})\times 
\mathbf{H}^{\frac{1}{2}}(\partial \mathcal{O})}=0\text{ \ for every }\mathbf{%
\tau }\in TH^{1/2}(\partial \mathcal{O})
\end{equation*}%
(and so $\left( \sigma (\mathbf{u}_{0})\mathbf{n}-p_{0}\mathbf{n}\right)
\cdot \mathbf{\tau }=0$ in the sense of distributions; see Remark 3.1 of 
\cite{agw}).

\item[(A.v)] The flow velocity component $\mathbf{u}_{0}=\mathbf{f}_{0}+%
\widetilde{\mathbf{f}}_{0}$, where $\mathbf{f}_{0}\in \mathbf{V}_{0}$ and $%
\widetilde{\mathbf{f}}_{0}\in \mathbf{H}^{1}(\mathcal{O})$ satisfies%
\footnote{%
The existence of an $\mathbf{H}^{1}(\mathcal{O})$-function $\widetilde{%
\mathbf{f}}_{0}$ with such a boundary trace on Lipschitz domain $\mathcal{O}$
is assured; see e.g., Theorem 3.33 of \cite{Mc}.}%
\begin{equation*}
\left. \widetilde{\mathbf{f}}_{0}\right\vert _{\partial \mathcal{O}}=%
\begin{cases}
0 & ~\text{ on }~S \\ 
w_{2}\mathbf{n} & ~\text{ on}~\Omega%
\end{cases}%
\end{equation*}%
\noindent (and so $\left. \mathbf{f}_{0}\right\vert _{\partial \mathcal{O}%
}\in TH^{1/2}(\partial \mathcal{O})$).
\end{enumerate}

\bigskip

With respect to the flow-structure dynamics (\ref{system1})--(\ref{IC_2}),
we established the following result in \cite{agw} under a standard \
assumption -- see p.529 of \cite{dV} and pp.102--103 of \cite{valli}) -- on
ambient (real-valued) vector field $\mathbf{U}$ of (\ref{system1}):

\medskip

\begin{theorem}
\label{well}(a) With given ambient vector field $\mathbf{U}\in \mathbf{V}%
_{0}\cap \mathbf{H}^{3}(\mathcal{O})$, the flow-structure operator $\mathcal{%
A}:D(\mathcal{A})\subset \mathcal{H}\rightarrow \mathcal{H}$, generates a $%
C_{0}$-semigroup on $\mathcal{H}$. Accordingly, the solution of (\ref%
{system1})--(\ref{IC_2}) -- with initial data $[p_{0},\mathbf{u}%
_{0},w_{0},w_{1}]\in $ $\mathcal{H}$ -- may be given by 
\begin{equation}
\begin{bmatrix}
p(t) \\ 
\mathbf{u}(t) \\ 
w(t) \\ 
w_{t}(t)%
\end{bmatrix}%
=e^{\mathcal{A}t}%
\begin{bmatrix}
p_{0} \\ 
\mathbf{u}_{0} \\ 
w_{0} \\ 
w_{1}%
\end{bmatrix}%
\in C[(0,\infty );\mathcal{H}].  \label{w_1}
\end{equation}%
Moreover, the flow-structure semigroup obeys the following estimate:

\begin{equation}
\left\Vert e^{\mathcal{A}t}\right\Vert _{\mathcal{L}(\mathcal{H})}\leq \exp
\left( \frac{t}{2}\left\Vert \func{div}(\mathbf{U})\right\Vert _{\infty
}\right) \text{, for every }t>0.  \label{w_2}
\end{equation}
\end{theorem}

\medskip

\section{Statement of Main Result and Literature}

\smallskip

\subsection{The Main Result}

Concerning the ambient vector field $\mathbf{U}$, we will assume throughout
that

\begin{equation}
\mathbf{U}\in \mathbf{V}_{0}\cap \mathbf{H}^{3}(\mathcal{O})\text{, and
further satisfies }\func{div}(\mathbf{U})=0\text{ in }\mathcal{O}.
\label{extra}
\end{equation}

\medskip

It will be shown below that the Null Space of $\mathcal{A}:D(\mathcal{A}%
)\subset \mathcal{H}\rightarrow \mathcal{H}$, denoted $Null(\mathcal{A)}$,
is nonempty (but only \emph{one dimensional}, with explicitly identified
eigenfunction; see Lemma \ref{SS}(a)). Moreover, one can cleanly
characterize the orthogonal complement as follows:%
\begin{equation}
\lbrack Null(\mathcal{A)}]^{\bot }=\mathcal{\{}[p_{0},\mathbf{u}%
_{0},w_{1},w_{2}]\in \mathcal{H}:\int\limits_{\mathcal{O}}p_{0}d\mathcal{O}%
+\int\limits_{\Omega }w_{1}d\Omega =0\mathcal{\}}  \label{N22}
\end{equation}%
(See Lemma \ref{SS}(b)). Accordingly, so as to avoid the possibility of
\textquotedblleft steady states\textquotedblright , our main result of
uniform decay for the compressible flow-structure PDE model (\ref{system1}%
)--(\ref{IC_2}) is applicable for initial data drawn from $[Null(\mathcal{A)}%
]^{\bot }$.

\begin{theorem}
\label{uniform} Assume that geometries $\left[ \mathcal{O},\Omega \right] $
are in either class $\mathsf{(G.1)}$ or $\mathsf{(G.2)}$. In addition,
assume that ambient vector field $\mathbf{U}$ satisfies (\ref{extra}). Then
solutions of (\ref{system1})--(\ref{IC_2}) decay uniformly in time for all
initial data in $[Null(\mathcal{A)}]^{\bot }$. That is, there exist
constants $M\geq 1$ and $\delta >0$ such that for all $t\geq 0$, 
\begin{equation*}
\left\Vert p(t),u(t),w(t),w_{t}(t)\right\Vert _{\mathcal{H}}\leq Me^{-\delta
t}\left\Vert [p_{0},\mathbf{u}_{0},w_{1},w_{2}]\right\Vert _{\mathcal{H}}%
\text{, \ \ for all }[p_{0},\mathbf{u}_{0},w_{1},w_{2}]\in \lbrack Null(%
\mathcal{A)}]^{\bot }.
\end{equation*}
\end{theorem}

\medskip

\begin{remark}
In \cite{Igor-note}, I. Chueshov stated that he was most interested in the
situation when $\mathcal{O}$ is a tube domain along the $x_{3}$-axis and
ambient vector $\mathbf{U}=[0,0,U]$, where constant (speed) $U>0$%
,\textquotedblleft ...the case which is important from the point of view of
aeroelasticity.\textquotedblright\ (See \cite{dowell1}, and also \cite%
{supersonic}, where $\mathbf{U}$ in the latter reference assumes this
particular form.) With this physical situation in mind, one should then take
ambient field $\mathbf{U}$ in (\ref{system1})--(\ref{IC_2}) to be 
\begin{equation*}
\mathbf{U}=U\mathbf{\phi }_{\epsilon }\text{,}
\end{equation*}%
where divergence free test function $\mathbf{\phi }_{\epsilon }\in \mathbf{C}%
_{0}^{\infty }(\mathcal{O})$ satisfies $\mathbf{\phi }_{\epsilon }(\mathbf{x}%
)=[0,0,1]$ \ if $dist(\mathbf{x},\partial \mathcal{O})>\epsilon $, and dies
off as $dist(\mathbf{x},\partial \mathcal{O})$ goes to zero.
\end{remark}

\medskip

\begin{remark}
The imposition in Theorem \ref{uniform} on the geometries $\left[ \mathcal{O}%
,\Omega \right] $ are in large part to eventually allow our properly
applying the higher regularity results in \cite{dauge}, results which
pertain to solutions of inhomogeneous Stokes equations on domains with
corners (see (\ref{5.24}) and (\ref{5.26}) below).
\end{remark}

\medskip

\medskip

\subsection{Literature}

With respect to the main issue of uniform decay and associated spectral
analysis for the coupled compressible flow-structure PDE system (\ref%
{system1})--(\ref{IC_2}), much of the relevant literature concerns stability
properties of \emph{uncoupled} compressible flows on \emph{unbounded }%
domains. Such works include \cite{deckel},\cite{kagei},\cite{kob},\cite%
{spectral}, \cite{decay} and references therein. For example, in \cite%
{spectral} and \cite{decay}, the three dimensional geometry in which their
(uncoupled) compressible flow PDE model evolves is of the form $D\times
(-\infty ,\infty )$, where $D$ is a bounded and connected domain: in these
papers, the authors are able to apply the Fourier transform in the third
variable, by way of obtaining their spectral and uniform decay results,
provided the associated Reynolds and Mach numbers are small enough. Since
our flow geometry $\mathcal{O}$ is generally assumed to be bounded and
convex, and no size requirements are imposed on our ambient vector field $%
\mathbf{U}$, such an approach cannot be availed of in the present case, at
least to resolve the flow component of the dynamics.

\smallskip

The one reference which pertains to exponential stability of the
flow-structure system (\ref{system1})--(\ref{IC_2}), in the special case $%
\mathbf{U}\equiv \mathbf{0}$, is the work \cite{Chu2013-comp} (This paper
also resolves wellposedness of [structurally] nonlinear versions of (\ref%
{system1})--(\ref{IC_2}) in the case of zero ambient state.) However, the
time dependent multiplier method which is adopted in \cite{Chu2013-comp}, by
way of establishing exponential decay, is \emph{inapplicable} in the present
case of generally nonzero vector field $\mathbf{U}$ which satisfies (\ref%
{extra}).

In particular, in \cite{Chu2013-comp} the following variable $N_{0}(p,w)$ is
constructed, with respect to the flow-structure system (\ref{system1})--(\ref%
{IC_2}), and with initial data taken from $[Null(\mathcal{A)}]^{\bot }$, as
characterized in (\ref{N22}):%
\begin{equation*}
N_{0}(p(t),w_{t}(t))=\nabla \psi \text{ }\Longleftrightarrow \left\{ 
\begin{array}{l}
\Delta \psi =-p\text{ \ in }\mathcal{O} \\ 
\frac{\partial }{\partial \mathbf{n}}\psi =0\text{ \ on }S\text{; \ }\frac{%
\partial }{\partial \mathbf{n}}\psi =w\text{ \ on }\Omega .%
\end{array}%
\right.
\end{equation*}%
(We note that from the invariance result in Proposition \ref{invariant}(b)
below, and the characterization (\ref{N22}), the data $\left\{
-p(t),w(t)\right\} $ for this BVP satisfies the compatibility condition
necessary for $\psi $ to be well-defined.) This variable is used in the
context of a multiplier method on page 660 of \cite{Chu2013-comp}, by way of
obtaining the needed $L^{2}(\mathcal{O})$ estimate for variable $p(t)$,
pointwise in time. However, this method uses very critically the fact that
ambient vector\ field $\mathbf{U}\equiv \mathbf{0}$, so as to exploit the
consequent relation $p_{t}=\func{div}(\mathbf{u})$, which comes from the
pressure state equation in (\ref{system1})--(\ref{IC_2}) with $\mathbf{U}%
\equiv \mathbf{0}$ therein. (See the middle of page 660 of \cite%
{Chu2013-comp}.) For generally nonzero $\mathbf{U}$, this argument in (\ref%
{system1})--(\ref{IC_2}) is unavailing, because of the presence of term $%
\mathbf{U}\cdot \nabla p$ in the pressure equation of (\ref{system1}).)

\medskip

Considering the aforesaid complication, we choose to operate in the
\textquotedblleft frequency domain\textquotedblright , instead of the time
domain: Namely, for all $\beta \in \mathbb{R}$, we obtain a uniform estimate
on the operator norm $\left( i\beta -\mathcal{A}\right) ^{-1}$ when
restricted to $[Null(\mathcal{A)}]^{\bot }$ as characterized in (\ref{N22})
(recall that $Null(\mathcal{A)}$ as described in Lemma \ref{SS}(a) below is 
\emph{one dimensional}). Such a uniform bound on the resolvent, as it acts
on the imaginary axis, ultimately allows for the wellknown resolvent
characterization for exponential stability of bounded semigroups; see
Theorem \ref{stab-un} in the Appendix, as well as \cite{huang} and \cite%
{pruss}. (Such a static methodology was previouslyused in \cite{ALT}, \cite%
{george1} and \cite{george2} in the context of determining uniform decay
rates for other compressible fluid-structure interactions.) The work to
estimate said resolvent is undertaken in Section \ref{freq} below.

\medskip

One definite benefit of the the time independent methodology which we employ
here, is that it ultimately allows for an adequate treatment of the term $%
\mathbf{U}\cdot \nabla p$, as it appears in the pressure equation of (\ref%
{system1}). This is accomplished in Section \ref{freq}, by means of an
appropriate decomposition of the \textquotedblleft zero mean
average\textquotedblright\ component of the pressure solution component, in
combination with known higher regularity results for inhomogeneous Stokes
systems on nonsmooth domains; see \cite{dauge}.

\medskip

\section{Some Supporting Results}

\subsection{A Basic Energy Equality in the Frequency Domain}

We begin by providing the following flow-structure (static) dissipation
relation.

\smallskip

\begin{proposition}
\label{energy}Let 
\begin{equation}
\Phi \equiv 
\begin{bmatrix}
p_{0} \\ 
\mathbf{u}_{0} \\ 
w_{1} \\ 
w_{2}%
\end{bmatrix}%
\in D(\mathcal{A)},  \label{short}
\end{equation}%
with $D(\mathcal{A)}$ being as described in (A.i)-(A.v) above. Then if the
ambient vector field $\mathbf{U}$ satisfies (\ref{extra}), we have 
\begin{equation}
\func{Re}\left( \mathcal{A}\Phi ,\Phi \right) _{\mathcal{H}}=-\left( \sigma (%
\mathbf{u}_{0}),\epsilon (\mathbf{u}_{0})\right) _{\mathcal{O}}-\eta
\left\Vert \mathbf{u}_{0}\right\Vert _{\mathcal{O}}^{2}.  \label{dissi}
\end{equation}
\end{proposition}

\medskip

\textbf{Proof of Proposition \ref{energy}:} Using the definition of $%
\mathcal{A}:D(\mathcal{A})\subset \mathcal{H}\rightarrow \mathcal{H}$ in (%
\ref{AAA}) and its domain $D(\mathcal{A)}$, we have 
\begin{equation}
\begin{array}{l}
\left( \mathcal{A}\Phi ,\Phi \right) _{\mathcal{H}}= \\ 
\text{ \ \ \ \ \ }\left( -\mathbf{U}\cdot \nabla p_{0}-div~\mathbf{u}%
_{0},p_{0}\right) _{\mathcal{O}}+\left( -\nabla p_{0}+\func{div}\sigma (%
\mathbf{u}_{0})-\eta \mathbf{u}_{0}-\mathbf{U}\cdot \nabla \mathbf{u}_{0},%
\mathbf{u}_{0}\right) _{\mathcal{O}} \\ 
\\ 
\text{ \ \ \ \ \ }+\left( \Delta w_{2},\Delta w_{1}\right) _{\Omega }+\left(
-\Delta ^{2}w_{1}-\left[ 2\nu \partial _{x_{3}}(\mathbf{u}_{0})_{3}+\lambda 
\text{div}(\mathbf{u}_{0})-p_{0}\right] _{\Omega },w_{2}\right) _{\Omega }.%
\end{array}
\label{int}
\end{equation}

\bigskip

Upon an integration by parts -- which uses the fact that ambient vector
field $\mathbf{U}\in \mathbf{V}_{0}$-- we then have

\begin{equation}
\begin{array}{l}
\left( \mathcal{A}\Phi ,\Phi \right) _{\mathcal{H}}=\left\langle \sigma (%
\mathbf{u}_{0})\cdot \mathbf{n-p}_{0}\mathbf{n},\mathbf{u}_{0}\right\rangle
_{\partial \mathcal{O}}-\left( \left[ 2\nu \partial _{x_{3}}(\mathbf{u}%
_{0})_{3}+\lambda \text{div}(\mathbf{u}_{0})-p_{0}\right] _{\Omega
},w_{2}\right) _{\Omega } \\ 
\\ 
-\left( \sigma (\mathbf{u}_{0}),\epsilon (\mathbf{u}_{0})\right) _{\mathcal{O%
}}-\eta \left\Vert \mathbf{u}_{0}\right\Vert _{\mathcal{O}}^{2}+\frac{1}{2}%
\left( \func{div}(\mathbf{U})p_{0},p_{0}\right) _{\mathcal{O}}+\frac{1}{2}%
\left( \func{div}(\mathbf{U})\mathbf{u}_{0},\mathbf{u}_{0}\right) _{\mathcal{%
O}} \\ 
\\ 
-i\func{Im}\left( \mathbf{U}\cdot \nabla p_{0},p_{0}\right) _{\mathcal{O}}-i%
\func{Im}\left( \mathbf{U}\cdot \nabla \mathbf{u}_{0},\mathbf{u}_{0}\right)
_{\mathcal{O}}-2i\func{Im}\left( \func{div}(\mathbf{u}_{0}),p_{0}\right) _{%
\mathcal{O}}-2i\func{Im}\left( \Delta w_{1},\Delta w_{2}\right) _{\Omega }.%
\end{array}
\label{int2}
\end{equation}%
Now since $\Phi =[p_{0},\mathbf{u}_{0},w_{1},w_{2}]\in D(\mathcal{A)}$, then
from (A.v), we have $\mathbf{u}_{0}=\mathbf{f}_{0}+\widetilde{\mathbf{f}}%
_{0} $, where $\mathbf{f}_{0}\in \mathbf{V}_{0}$ and $\widetilde{\mathbf{f}}%
_{0}\in \mathbf{H}^{1}(\mathcal{O})$ satisfies%
\begin{equation}
\left. \widetilde{\mathbf{f}}_{0}\right\vert _{\partial \mathcal{O}}=%
\begin{cases}
0 & ~\text{ on }~S \\ 
w_{2}\mathbf{n} & ~\text{ on}~\Omega .%
\end{cases}
\label{int3}
\end{equation}%
Since component $\left. \mathbf{f}_{0}\right\vert _{\partial \mathcal{O}}\in
TH^{1/2}(\partial \mathcal{O})$, then using (A.iv), (\ref{int3}), and the
fact that $\left. \mathbf{n}\right\vert _{\Omega }=(0,0,1)$, we have that
for the first two terms on RHS of (\ref{int2}), 
\begin{eqnarray*}
&&\left\langle \sigma (\mathbf{u}_{0})\cdot \mathbf{n-p}_{0}\mathbf{n},%
\mathbf{u}_{0}\right\rangle _{\partial \mathcal{O}} \\
&=&\left\langle \sigma (\mathbf{u}_{0})\cdot \mathbf{n-p}_{0}\mathbf{n},%
\widetilde{\mathbf{f}}_{0}\right\rangle _{\partial \mathcal{O}}-\left( \left[
2\nu \partial _{x_{3}}(\mathbf{u}_{0})_{3}+\lambda \text{div}(\mathbf{u}%
_{0})-p_{0}\right] _{\Omega },w_{2}\right) _{\Omega } \\
&=&0.
\end{eqnarray*}

Applying this relation and using the assumption that $\mathbf{U}$ is
solenoidal, (\ref{int2}) then becomes%
\begin{eqnarray}
&&\left( \mathcal{A}\Phi ,\Phi \right) _{\mathcal{H}}=-\left( \sigma (%
\mathbf{u}_{0}),\epsilon (\mathbf{u}_{0})\right) _{\mathcal{O}}-\eta
\left\Vert \mathbf{u}_{0}\right\Vert _{\mathcal{O}}^{2}  \notag \\
&&-i\func{Im}\left( \mathbf{U}\cdot \nabla p_{0},p_{0}\right) _{\mathcal{O}%
}-i\func{Im}\left( \mathbf{U}\cdot \nabla \mathbf{u}_{0},\mathbf{u}%
_{0}\right) _{\mathcal{O}}-2i\func{Im}\left( \func{div}(\mathbf{u}%
_{0}),p_{0}\right) _{\mathcal{O}}-2i\func{Im}\left( \Delta w_{1},\Delta
w_{2}\right) _{\Omega }.  \label{int4}
\end{eqnarray}

This establishes (\ref{dissi}). \ \ \ $\square $

\medskip

\subsection{A Spectral Analysis on $i\mathbb{R}$}

\bigskip

\subsubsection{Concerning the point $\protect\lambda =0$}

\begin{lemma}
\label{SS}Under the assumption that ambient vector field $\mathbf{U}$
satisfies (\ref{extra}), one has: (a) The subspace $Null(\mathcal{A)}\subset 
\mathcal{H}$ of the flow-structure generator $\mathcal{A}:D(\mathcal{A}%
)\subset \mathcal{H}\rightarrow \mathcal{H}$ is one dimensional. In
particular,%
\begin{equation}
Null(\mathcal{A)}=Span\left\{ \left[ 
\begin{array}{c}
1 \\ 
0 \\ 
{{\mathring{A}}^{-1}(1)} \\ 
0%
\end{array}%
\right] \right\} ,  \label{N1}
\end{equation}%
where $\mathring{A}:L^{2}(\Omega )\rightarrow L^{2}(\Omega )$ is the
following elliptic operator:%
\begin{equation}
\mathring{A}\varpi =\Delta ^{2}\varpi \text{, with }D(\mathring{A})=\{w\in
H_{0}^{2}(\Omega ):\Delta ^{2}w\in L^{2}(\Omega )\}.\text{\ \ }
\label{angst}
\end{equation}

(b) The orthogonal complement of $Null(\mathcal{A)}$ admits of the following
characterization: 
\begin{equation}
\lbrack Null(\mathcal{A)}]^{\bot }=\mathcal{\{}[p_{0},\mathbf{u}%
_{0},w_{1},w_{2}]\in \mathcal{H}:\int\limits_{\mathcal{O}}p_{0}d\mathcal{O}%
+\int\limits_{\Omega }w_{1}d\Omega =0\mathcal{\}}.  \label{N2}
\end{equation}
\end{lemma}

\mathstrut

\textbf{Proof:} Let $\Phi =[p_{0},\mathbf{u}_{0},w_{1},w_{2}]\in Null(%
\mathcal{A)}$, where $D(\mathcal{A)}$ is as given in (A.i)-(A.v) above; viz.,%
\begin{equation}
\mathcal{A}\Phi =0\text{.}  \label{nu_0}
\end{equation}%
This relation and (\ref{dissi}) (and Korn's Inequality) then give%
\begin{equation}
\mathbf{u}_{0}=0\text{.}  \label{nu_1}
\end{equation}%
In turn, from (\ref{nu_0}) and the second row of the matrix $\mathcal{A}$ as
given in (\ref{AAA}), we have $\nabla p_{0}=0$; so 
\begin{equation}
p_{0}=C_{0}\text{, \ where }C_{0}\equiv \text{ constant.}  \label{nu_1.5}
\end{equation}%
Moreover, from (A.v) and (\ref{nu_1}%
\begin{equation}
w_{2}=0\text{.}  \label{nu_2}
\end{equation}%
Lastly, from the fourth row of the matrix $\mathcal{A}$ and (\ref{nu_1.5})
we see that structural displacement variable $w_{1}$ satisfies the boundary
value problem%
\begin{equation*}
\Delta ^{2}w_{1}=C_{0}\text{ in }\Omega \text{, \ }\left. w_{1}\right\vert
_{\partial \Omega }=\left. \frac{\partial w_{1}}{\partial \nu }\right\vert
_{\partial \Omega }=0\text{ \ on }\partial \Omega \text{,}
\end{equation*}%
and so necessarily%
\begin{equation}
w_{1}={{\mathring{A}}^{-1}(C_{0}),}  \label{nu_3}
\end{equation}%
where $\mathring{A}^{-1}:L^{2}(\Omega )\rightarrow L^{2}(\Omega )$ is the
operator of the positive definite, self-adjoint operator defined in (\ref%
{angst}). The collection of (\ref{nu_1}), (\ref{nu_1.5}), (\ref{nu_2}) and (%
\ref{nu_3}) establishes (a). Given the definition of the $\mathcal{H}$-inner
product in (\ref{innerp}), as well as the definition (\ref{angst}), the
relation in (b) is immediate. \ \ \ $\square $

\medskip

\subsubsection{Concerning Discrete Spectra on $i%
\mathbb{R}
\diagdown \{0\}$}

\medskip

\begin{lemma}
\label{point_2} Under the assumption that ambient vector field $\mathbf{U}$
satisfies (\ref{extra}): For given $\beta \in 
\mathbb{R}
\diagdown \{0\}$, \newline
$Null(i\beta I-\mathcal{A})=\left\{ 0\right\} $.
\end{lemma}

We consider the equation%
\begin{equation}
\mathcal{A}%
\begin{bmatrix}
p_{0} \\ 
\mathbf{u}_{0} \\ 
w_{1} \\ 
w_{2}%
\end{bmatrix}%
=i\beta 
\begin{bmatrix}
p_{0} \\ 
\mathbf{u}_{0} \\ 
w_{1} \\ 
w_{2}%
\end{bmatrix}%
,\text{ \ \ where }\beta \neq 0\text{ and }%
\begin{bmatrix}
p_{0} \\ 
\mathbf{u}_{0} \\ 
w_{1} \\ 
w_{2}%
\end{bmatrix}%
\in D(\mathcal{A})\text{.}  \label{point}
\end{equation}

This abstract equation may equivalently be written as 
\begin{eqnarray}
&&\left\{ 
\begin{array}{l}
-\mathbf{U}\cdot \nabla p_{0}-div~\mathbf{u}_{0}=i\beta p_{0}\text{ \ in }~%
\mathcal{O} \\ 
-\nabla p_{0}+div~\sigma (\mathbf{u}_{0})-\eta \mathbf{u}_{0}-\mathbf{U}%
\cdot \nabla \mathbf{u}_{0}=i\beta \mathbf{u}_{0}~\text{ \ in }~\mathcal{O}
\\ 
w_{2}=i\beta w_{1}\text{ \ on }~\Omega \\ 
-\Delta ^{2}w_{1}-\left[ 2\nu \partial _{x_{3}}(\mathbf{u}_{0})_{3}+\lambda 
\text{div}(\mathbf{u}_{0})-p_{0}\right] _{\Omega }=i\beta w_{2}~\text{ on }%
~\Omega%
\end{array}%
\right.  \label{sys2} \\
&&  \notag \\
&&\left\{ 
\begin{array}{l}
(\sigma (\mathbf{u}_{0})\mathbf{n}-p_{0}\mathbf{n})\cdot \boldsymbol{\tau }=0%
\text{ \ on }~\partial \mathcal{O} \\ 
\mathbf{u}_{0}\cdot \mathbf{n}=0~\text{ on }~S \\ 
\mathbf{u}_{0}\cdot \mathbf{n}=w_{2}~\text{ on }~\Omega \\ 
w_{1}=\frac{\partial w_{1}}{\partial \nu }=0~\text{ on }~\partial \Omega%
\end{array}%
\right.  \label{sys3}
\end{eqnarray}%
Taking the $\mathcal{H}$-inner product of both sides of (\ref{point}) with
respect to $\Phi \equiv \lbrack p_{0},\mathbf{u}_{0},w_{1},w_{2}]$, and
subsequently applying the energy relation in Proposition \ref{energy} (along
with Korn's Inequality), we then infer

\begin{equation}
\mathbf{u}_{0}=0.  \label{b_0}
\end{equation}%
In turn, we get from the flow equation in (\ref{sys2}) that 
\begin{equation*}
\nabla p_{0}=0\Rightarrow p_{0}=\text{constant.}
\end{equation*}%
If we take this into account in the pressure equation in (\ref{sys2}), then
from (\ref{b_0}) we have 
\begin{equation}
p_{0}=0,~~~\text{as}~~~\beta \neq 0.  \label{b_1}
\end{equation}%
Also, using the normal component boundary conditions in (\ref{sys3}), we see
that 
\begin{equation}
w_{2}=\mathbf{u}_{0}\cdot n=0.  \label{b_2}
\end{equation}%
Applying this equality in turn to the structural equation in (\ref{sys2}),
we have finally (as $w_{1}$ satisfies the essential B.C.'s in (\ref{sys3}))%
\begin{equation}
{{\mathring{A}}}w_{1}=0~~\text{ in}~~\Omega \Rightarrow w_{1}=0.  \label{b_3}
\end{equation}%
Collecting (\ref{b_0})-(\ref{b_3}) now gives the asserted statement. \ \ \ $%
\square $

\bigskip

\subsubsection{A Certain A Priori Estimate}

\medskip

We derive here the following result which will be invoked throughout.

\smallskip

\begin{lemma}
\label{med}Under the assumption that ambient vector field $\mathbf{U}$
satisfies (\ref{extra}): Suppose that for given parameter $\beta \in 
\mathbb{R}
$ and given $\Phi ^{\ast }=[p_{0}^{\ast },\mathbf{u}_{0}^{\ast },w_{1}^{\ast
},w_{2}^{\ast }]\in \mathcal{H}$, $\Phi =[p_{0},\mathbf{u}%
_{0},w_{1},w_{2}]\in D(\mathcal{A})$ is a solution of 
\begin{equation}
(i\beta I-\mathcal{A})\Phi =\Phi ^{\ast }.  \label{lam}
\end{equation}%
Therewith, denote $L^{2}$-function $q_{0}$ to be \textquotedblleft the zero
average\textquotedblright\ component of pressure variable $q_{0}$; viz., $%
q_{0}$ is from the decomposition 
\begin{equation}
p_{0}=q_{0}+c_{0}\text{,}  \label{q}
\end{equation}

where 
\begin{equation}
q_{0}\text{ satisfies }\int\limits_{\mathcal{O}}q_{0}d\mathcal{O}=0\text{, \
and \ }c_{0}=\text{constant.}  \label{c0}
\end{equation}%
Then one
\end{lemma}

\begin{lemma}
has the estimate 
\begin{equation}
\left\Vert q_{0}\right\Vert _{\mathcal{O}}+\left\Vert \sigma (\mathbf{u}_{0})%
\mathbf{n}-q_{0}\mathbf{n}\right\Vert _{\mathbf{H}^{-\frac{1}{2}}(\partial 
\mathcal{O})}\leq C\left( \sqrt{\left\Vert \Phi \right\Vert _{\mathcal{H}%
}\left\Vert \Phi ^{\ast }\right\Vert _{\mathcal{H}}}+\left\Vert \Phi ^{\ast
}\right\Vert _{\mathcal{H}}+\left\vert \beta \right\vert \left\Vert \mathbf{u%
}_{0}\right\Vert _{\mathcal{O}}\right) .  \label{med2}
\end{equation}
\end{lemma}

\medskip

\textbf{Proof:} To start, we take the $\mathcal{H}$-inner product of both
sides of (\ref{lam}) with respect to $\Phi $, and subsequently appeal to
Proposition \ref{energy} -- with assumption (\ref{extra}) in play -- so as
to have%
\begin{equation}
\left( \sigma (\mathbf{u}_{0}),\epsilon (\mathbf{u}_{0})\right) _{\mathcal{O}%
}+\eta \left\Vert \mathbf{u}_{0}\right\Vert _{\mathcal{O}}^{2}\leq
\left\vert \left( \Phi ^{\ast },\Phi \right) _{\mathcal{H}}\right\vert .
\label{lam2}
\end{equation}%
In turn, from (\ref{lam}) and (\ref{q}), $\left\{ \mathbf{u}%
_{0},q_{0}\right\} $ satisfies the boundary problem 
\begin{equation}
\left\{ 
\begin{array}{l}
-div~\sigma (\mathbf{u}_{0})+\eta \mathbf{u}_{0}+\nabla q_{0}=\mathbf{u}%
_{0}^{\ast }-i\beta \mathbf{u}_{0}-\mathbf{U}\cdot \nabla \mathbf{u}_{0}~%
\text{ in }~\mathcal{O} \\ 
\func{div}(\mathbf{u}_{0})=\func{div}(\mathbf{u}_{0})\text{ \ in }\mathcal{O}
\\ 
\left. \mathbf{u}_{0}\right\vert _{\partial \mathcal{O}}=\left. \mathbf{u}%
_{0}\right\vert _{\partial \mathcal{O}}\text{ on }\partial \mathcal{O}%
\end{array}%
\right.  \label{lam3}
\end{equation}

(note that because of (\ref{lam2}), we have control of the data of this
BVP). Accordingly, we can appeal to classic (incompressible) Stokes theory
(see Theorem \ref{static} of the Appendix), so as to have the estimate 
\begin{eqnarray*}
\left\Vert q_{0}\right\Vert _{\mathcal{O}} &\leq &C\left[ \left\Vert \mathbf{%
u}_{0}\right\Vert _{\mathbf{H}^{1}(\mathcal{O})}+\left\Vert \mathbf{u}%
_{0}^{\ast }\right\Vert _{\mathcal{O}}+\left\vert \beta \right\vert
\left\Vert \mathbf{u}_{0}\right\Vert _{\mathcal{O}}\right. \\
&&+\left. \left\Vert \mathbf{U}\cdot \nabla \mathbf{u}_{0}\right\Vert _{%
\mathcal{O}}+\left\Vert \func{div}(\mathbf{u}_{0})\right\Vert _{\mathcal{O}%
}+\left\Vert \left. \mathbf{u}_{0}\right\vert _{\partial \mathcal{O}%
}\right\Vert _{H^{\frac{1}{2}}(\partial \mathcal{O})}\right] .
\end{eqnarray*}%
An estimation of right hand side, by means of (\ref{lam2}), Korn's
Inequality, and the Sobolev Trace Theorem then gives

\begin{equation}
\left\Vert q_{0}\right\Vert _{\mathcal{O}}\leq C\left( \sqrt{\left\Vert \Phi
\right\Vert _{\mathcal{H}}\left\Vert \Phi ^{\ast }\right\Vert _{\mathcal{H}}}%
+\left\Vert \Phi ^{\ast }\right\Vert _{\mathcal{H}}+\left\vert \beta
\right\vert \left\Vert \mathbf{u}_{0}\right\Vert _{\mathcal{O}}\right) .
\label{lam4}
\end{equation}

\smallskip

Moreover, a standard energy method with respect to the flow equation in (\ref%
{lam3}) -- which uses the surjectivity and boundedness of the Sobolev Trace
Map for $\mathbf{H}^{1}$-functions on Lipschitz domains -- see e.g., Theorem
3.38 of \cite{Mc} -- eventually yields the following boundary trace estimate
for solution components $\left\{ \mathbf{u}_{0},q_{0}\right\} $: 
\begin{equation*}
\left\Vert \sigma (\mathbf{u}_{0})\mathbf{n}-q_{0}\mathbf{n}\right\Vert _{%
\mathbf{H}^{-\frac{1}{2}}(\partial \mathcal{O})}\leq C\left[ \left\Vert 
\mathbf{u}_{0}\right\Vert _{\mathbf{H}^{1}(\mathcal{O})}+\left\Vert
q_{0}\right\Vert _{\mathcal{O}}+\left\Vert div~\sigma (\mathbf{u}%
_{0})-\nabla q_{0}\right\Vert _{\mathcal{O}}\right]
\end{equation*}%
\begin{equation*}
=C\left[ \left\Vert \mathbf{u}_{0}\right\Vert _{\mathbf{H}^{1}(\mathcal{O}%
)}+\left\Vert q_{0}\right\Vert _{\mathcal{O}}+\left\Vert \mathbf{u}%
_{0}^{\ast }-~i\beta \mathbf{u}_{0}-\eta \mathbf{u}_{0}-\mathbf{U}\cdot
\nabla \mathbf{u}_{0}\right\Vert _{\mathcal{O}}\right] .
\end{equation*}%
Estimating this right hand side with (\ref{lam2}), (\ref{lam4}), and Korn's
Inequality, we now have%
\begin{equation}
\left\Vert \sigma (\mathbf{u}_{0})\mathbf{n}-q_{0}\mathbf{n}\right\Vert _{%
\mathbf{H}^{-\frac{1}{2}}(\partial \mathcal{O})}\leq C\left( \sqrt{%
\left\Vert \Phi \right\Vert _{\mathcal{H}}\left\Vert \Phi ^{\ast
}\right\Vert _{\mathcal{H}}}+\left\Vert \Phi ^{\ast }\right\Vert _{\mathcal{H%
}}+\left\vert \beta \right\vert \left\Vert \mathbf{u}_{0}\right\Vert _{%
\mathcal{O}}\right) .  \label{lam5}
\end{equation}%
Adding estimates (\ref{lam4}) and (\ref{lam5}) now completes the proof. \ \
\ $\square $

\subsubsection{Concerning Continuous Spectra on $i%
\mathbb{R}
\diagdown \{0\}$}

\medskip

The driving result in this Section is the following:

\smallskip

\begin{lemma}
\label{approx}Under the assumption that ambient vector field $\mathbf{U}$
satisfies (\ref{extra}): There is no approximate spectrum of $\mathcal{A}:D(%
\mathcal{A})\subset \mathcal{H}\rightarrow \mathcal{H}$ on $i%
\mathbb{R}
\diagdown \{0\}$.
\end{lemma}

\smallskip

\textbf{Proof: }Suppose that $i\beta $ is in the approximate spectrum of $%
\mathcal{A}$. Then by definition -- see. e.g., \cite{friedman}, p. 127 --
there exists a sequence 
\begin{equation*}
\left\{ \Phi _{n}\right\} \equiv \left\{ 
\begin{bmatrix}
p_{n} \\ 
\mathbf{u}_{n} \\ 
w_{1,n} \\ 
w_{2,n}%
\end{bmatrix}%
\right\} \subset D(\mathcal{A})
\end{equation*}%
such that for every $n\in \mathbb{N}$, 
\begin{equation}
\left\Vert \Phi _{n}\right\Vert _{\mathcal{H}}=1\text{ }  \label{ap1}
\end{equation}%
and 
\begin{eqnarray}
\left\Vert \Phi _{n}^{\ast }\right\Vert _{\mathcal{H}} &<&\frac{1}{n}\text{,}
\notag \\
\text{where }\Phi _{n}^{\ast } &=&\left[ p_{n}^{\ast },\mathbf{u}_{n}^{\ast
},w_{1,n}^{\ast },w_{2,n}^{\ast }\right] \equiv (i\beta I-\mathcal{A)}\Phi
_{n}\in \mathcal{H}\text{.}  \label{ap2}
\end{eqnarray}%
With the above notation, we then have the following static flow-structure
system:

\begin{equation}
(i\beta I-\mathcal{A)}\Phi _{n}=\Phi _{n}^{\ast },  \label{ap3}
\end{equation}%
or equivalently, from the form of the matrix $\mathcal{A}$ in (\ref{AAA})

\begin{align}
& \left\{ 
\begin{array}{l}
i\beta p_{n}+\mathbf{U}\cdot \nabla p_{n}+div~\mathbf{u}_{n}=p_{n}^{\ast }~%
\text{ in }~\mathcal{O} \\ 
i\beta \mathbf{u}_{n}+\nabla p_{n}-div~\sigma (\mathbf{u}_{n})+\eta \mathbf{u%
}_{n}+\mathbf{U}\cdot \nabla \mathbf{u}_{n}=\mathbf{u}_{n}^{\ast }~\text{ in 
}~\mathcal{O} \\ 
i\beta w_{1,n}-w_{2,n}=w_{1,n}^{\ast }\text{ \ in }~\Omega \\ 
i\beta w_{2,n}+\Delta ^{2}w_{1,n}+\left[ 2\nu \partial _{x_{3}}(\mathbf{u}%
_{n})_{3}+\lambda \text{div}(\mathbf{u}_{n})-p_{n}\right] _{\Omega
}=w_{2,n}^{\ast }~\text{ on }~\Omega%
\end{array}%
\right.  \label{ap4} \\
&  \notag \\
& \left\{ 
\begin{array}{l}
(\sigma (\mathbf{u}_{n})\mathbf{n}-p_{n}\mathbf{n})\cdot \boldsymbol{\tau }=0%
\text{ \ on }~\partial \mathcal{O} \\ 
\mathbf{u}_{n}\cdot \mathbf{n}=0~\text{ on }~S \\ 
\mathbf{u}_{n}\cdot \mathbf{n}=w_{2,n}~\text{ on }~\Omega \\ 
w_{1,n}=\frac{\partial w_{1,n}}{\partial \nu }=0~\text{ on }~\partial \Omega
.%
\end{array}%
\right.  \label{ap5}
\end{align}

\medskip

Taking the $\mathcal{H}$-inner product of both sides of (\ref{ap3}) with
respect to $\Phi _{n}$, and subsequently appealing to Proposition \ref%
{energy}, we have then%
\begin{equation*}
\left( \sigma (\mathbf{u}_{n}),\epsilon (\mathbf{u}_{n})\right) _{\mathcal{O}%
}+\eta \left\Vert \mathbf{u}_{n}\right\Vert _{\mathcal{O}}^{2}=\func{Re}%
\left( \Phi _{n}^{\ast },\Phi _{n}\right) _{\mathcal{H}},
\end{equation*}%
whence we obtain from (\ref{ap1}) and (\ref{ap2}),

\begin{equation}
\left( \sigma (\mathbf{u}_{n}),\epsilon (\mathbf{u}_{n})\right) +\eta
\left\Vert \mathbf{u}_{n}\right\Vert ^{2}=\mathcal{O}\left( \frac{1}{n}%
\right) .  \label{ap6}
\end{equation}%
Subsequently, for each $n$, we invoke again the wellknown $L^{2}-$%
decomposition of each $p_{n}.$ That is, we set 
\begin{equation}
p_{n}=q_{n}+c_{n},\text{ \ }\forall n,  \label{ap7}
\end{equation}%
where 
\begin{equation}
L^{2}\text{-function }q_{n}\text{ satisfies }\int\limits_{\mathcal{O}}q_{n}d%
\mathcal{O}=0\text{, \ and \ }c_{n}=\text{constant.}  \label{ap8}
\end{equation}%
For each $n$, the estimate (\ref{med2}) in Lemma \ref{med} directly applies;
we combine it with (\ref{ap1}), (\ref{ap2}), and (\ref{ap6}), so as to have
for fixed $\beta \neq 0$%
\begin{equation}
\left\Vert q_{n}\right\Vert _{\mathcal{O}}+\left\Vert \sigma (\mathbf{u}_{n})%
\mathbf{n}-q_{n}\mathbf{n}\right\Vert _{\mathbf{H}^{-\frac{1}{2}}(\partial 
\mathcal{O})}=O\left( \frac{1}{\sqrt{n}}\right) .  \label{ap9}
\end{equation}%
In turn, we can apply the decomposition (\ref{ap7}) to the pressure equation
in (\ref{ap4}), so as to have%
\begin{equation}
i\beta c_{n}=p_{n}^{\ast }-\mathbf{U}\cdot \nabla q_{n}-\func{div}\mathbf{u}%
_{n}-i\beta q_{n}\text{ \ in }\mathcal{O}.  \label{ap10}
\end{equation}%
Multiplying both sides of this equation by (constant) $\overline{c_{n}}$,
and then integrating by parts, we have%
\begin{equation}
i\beta \lbrack meas(\mathcal{O})]\left\vert c_{n}\right\vert ^{2}=\overline{%
c_{n}}\left[ \int\limits_{\mathcal{O}}p_{n}^{\ast }d\mathcal{O}-\int\limits_{%
\mathcal{O}}\func{div}(\mathbf{u}_{n})d\mathcal{O}\right]  \label{ap10.2}
\end{equation}%
(where we have implicitly used $\mathbf{U}\cdot \mathbf{n}=0$ on $\partial 
\mathcal{O}$, $\func{div}(\mathbf{U})=0$ in $\mathcal{O}$, and $\int\limits_{%
\mathcal{O}}q_{n}d\mathcal{O=}0$). To estimate right hand side of (\ref%
{ap10.2}), we note initially that from (\ref{ap7}) that%
\begin{equation*}
meas(\mathcal{O})\left\vert c_{n}\right\vert ^{2}=\left\Vert
p_{n}\right\Vert _{\mathcal{O}}^{2}-\left\Vert q_{n}\right\Vert _{\mathcal{O}%
}^{2},
\end{equation*}%
whence we obtain from (\ref{ap1}) and (\ref{ap9}),%
\begin{equation*}
\left\vert c_{n}\right\vert ^{2}\leq C_{\mathcal{O}}\left( 1+\frac{1}{n}%
\right) .
\end{equation*}%
Subsequently applying this estimate to the right hand side of (\ref{ap10.2}%
), along with (\ref{ap2}), and (\ref{ap6}) (and Cauchy-Schwartz), we then
have, for fixed $\beta \neq 0$, that actually%
\begin{equation}
\left\vert c_{n}\right\vert ^{2}=O\left( \frac{1}{n}\right) .  \label{c_n}
\end{equation}%
Combining this estimate with (\ref{ap7}) and (\ref{ap9}), we get now%
\begin{equation}
\left\Vert p_{n}\right\Vert _{\mathcal{O}}=O\left( \frac{1}{\sqrt{n}}\right)
.  \label{ap11}
\end{equation}

\smallskip

Moreover, using the normal component boundary condition in (\ref{ap5}), the
Sobolev Trace Theorem and (\ref{ap6}), we subsequently obtain%
\begin{equation}
\left\Vert w_{2,n}\right\Vert _{\Omega }=\left\Vert \mathbf{u}_{n}\cdot 
\mathbf{n\mid }_{\Omega }\right\Vert _{\Omega }=O\left( \frac{1}{\sqrt{n}}%
\right) .  \label{ap13}
\end{equation}

\smallskip

Lastly, we multiply both sides of the mechanical equation in (\ref{ap4}) by $%
w_{1,n}$ and integrate by parts to have%
\begin{equation}
\left\Vert \Delta w_{1,n}\right\Vert _{\Omega }^{2}=-\left\langle \sigma (%
\mathbf{u}_{n})\mathbf{n}-p_{n}\mathbf{n,}\left( w_{1,n}\mathbf{n}\right)
_{ext}\right\rangle _{\partial \mathcal{O}}-i\beta \left(
w_{2,n},w_{1,n}\right) _{\Omega }+\left( w_{2,n}^{\ast }~,w_{1,n}\right)
_{\Omega }.  \label{ap14}
\end{equation}%
Here, $\mathbf{H}^{2}(\partial \mathcal{O})-$function $\left( w_{1,n}\mathbf{%
n}\right) _{ext}$ is given by%
\begin{equation}
\left( w_{1,n}\mathbf{n}\right) _{ext}(x)=%
\begin{cases}
0 & ~\text{ }x\in \text{ }~S \\ 
w_{1,n}\mathbf{n} & ~\text{ }x\in ~\Omega%
\end{cases}%
.  \label{ap15}
\end{equation}%
Now, estimating right hand side of (\ref{ap14}) by (\ref{ap1}), (\ref{ap2}),
(\ref{ap7}), (\ref{ap9}), (\ref{c_n}) and (\ref{ap13}), we arrive at%
\begin{equation}
\left\Vert \Delta w_{1,n}\right\Vert _{\Omega }=\vartheta \left( \frac{1}{%
\sqrt{n}}\right) .  \label{ap16}
\end{equation}%
To conclude the proof: combining the estimates (\ref{ap6}), (\ref{ap11}), (%
\ref{ap13}) and (\ref{ap16}), we obtain%
\begin{equation*}
\left\Vert \Phi _{n}\right\Vert _{\mathcal{H}}=\left\Vert 
\begin{bmatrix}
p_{n} \\ 
\mathbf{u}_{n} \\ 
w_{1,n} \\ 
w_{2,n}%
\end{bmatrix}%
\right\Vert _{\mathcal{H}}\rightarrow 0.
\end{equation*}%
But this convergence contradicts (\ref{ap1}). Hence, $\lambda =i\beta $ is
not in the approximate spectrum of $\mathcal{A}$ for given $\beta \neq 0$. \
\ \ $\square $

\bigskip

From Lemma \ref{approx} we have immediately (see e.g., Theorem 2.27, p. 128
of \cite{friedman}),

\smallskip

\begin{corollary}
\label{coro}Under the assumption that ambient vector field $\mathbf{U}$
satisfies (\ref{extra}), one has\newline
$\sigma _{c}(\mathcal{A})\cap \left( i%
\mathbb{R}
\diagdown \{0\}\right) =\varnothing $.
\end{corollary}

\medskip

\subsubsection{Concerning $\protect\sigma _{r}(\mathcal{A)}\cap i%
\mathbb{R}
$}

The possibility of residual spectrum of the flow-structure generator on the
imaginary axis is eliminated quickly, after considering the representation
of the adjoint operator $\mathcal{A}^{\ast }:D(\mathcal{A}^{\ast })\subset 
\mathcal{H}\rightarrow \mathcal{H}$, and Lemma \ref{point_2}. In fact, a
standard computation yields,

\medskip

\begin{proposition}
\label{adj}Under the assumption that ambient vector field $\mathbf{U}$
satisfies (\ref{extra}), then the Hilbert space adjoint $\mathcal{A}^{\ast
}:D(\mathcal{A}^{\ast })\subset \mathcal{H}\rightarrow \mathcal{H}$ of $%
\mathcal{A}$ in (\ref{AAA}) is given by%
\begin{equation}
\mathcal{A}^{\ast }=\left[ 
\begin{array}{cccc}
\mathbf{U}\cdot \nabla (\cdot ) & \func{div}(\cdot ) & 0 & 0 \\ 
\nabla (\cdot ) & \func{div}\sigma (\cdot )-\eta (\cdot )+\mathbf{U}\cdot
\nabla (\cdot ) & 0 & 0 \\ 
0 & 0 & 0 & -I \\ 
-\left. \cdot \right\vert _{\Omega } & -[2\nu \partial _{x_{3}}(\cdot
)_{3}+\lambda \func{div}(\cdot )]_{\Omega } & \Delta ^{2} & 0%
\end{array}%
\right] .  \label{adjoint}
\end{equation}

Here, the domain $D(\mathcal{A}^{\ast })$ is given as 
\begin{equation*}
D(\mathcal{A}^{\ast })=\{(p_{0},\mathbf{u}_{0},w_{1},w_{2})\in L^{2}(%
\mathcal{O})\times \mathbf{H}^{1}(\mathcal{O})\times H_{0}^{2}(\Omega
)\times H_{0}^{2}(\Omega )~:~(A^{\ast }.i)\text{--}(A^{\ast }.v)~~\text{hold
below}\},
\end{equation*}%
where

(A$^{\ast }$.i) $\mathbf{U}\cdot \nabla p_{0}\in L^{2}(\mathcal{O})$

(A$^{\ast }$.ii) $\text{div}~\sigma (\mathbf{u}_{0})+\nabla p_{0}\in L^{2}(%
\mathcal{O})$

(A$^{\ast }$.iii) $\Delta ^{2}w_{0}-\left[ 2\nu \partial _{x_{3}}(\mathbf{u}%
_{0})_{3}+\lambda \text{div}(\mathbf{u}_{0})\right] _{\Omega }-\left.
p_{0}\right\vert _{\Omega }\in L^{2}(\Omega )$

(A$^{\ast }$.iv) $\left( \sigma (\mathbf{u}_{0})\mathbf{n}-p_{0}\mathbf{n}%
\right) \bot ~TH^{1/2}(\partial \mathcal{O})$. That is, 
\begin{equation*}
\left\langle \sigma (\mathbf{u}_{0})\mathbf{n}-p_{0}\mathbf{n},\mathbf{\tau }%
\right\rangle _{\mathbf{H}^{-\frac{1}{2}}(\partial \mathcal{O})\times 
\mathbf{H}^{\frac{1}{2}}(\partial \mathcal{O})}=0\text{ \ for all }\mathbf{%
\tau }\in TH^{1/2}(\partial \mathcal{O}).
\end{equation*}

(A$^{\ast }$.v) $\mathbf{u}_{0}=\mathbf{f}_{0}+\widetilde{\mathbf{f}}_{0}$,
where $\mathbf{f}_{0}\in \mathbf{V}_{0}$ and $\widetilde{\mathbf{f}}_{0}\in 
\mathbf{H}^{1}(\mathcal{O})$ satisfies 
\begin{equation*}
\left. \widetilde{\mathbf{f}}_{0}\right\vert _{\partial \mathcal{O}}=%
\begin{cases}
0 & ~\text{ on }~S \\ 
w_{2}\mathbf{n} & ~\text{ on}~\Omega%
\end{cases}%
\end{equation*}%
\noindent (and so $\left. \mathbf{f}_{0}\right\vert _{\partial \mathcal{O}%
}\in TH^{1/2}(\partial \mathcal{O})$).
\end{proposition}

\bigskip

\begin{lemma}
\label{adj_spec}Under the assumption that ambient vector field $\mathbf{U}$
satisfies (\ref{extra}), then:

(a)%
\begin{equation}
Null(\mathcal{A}^{\ast })=Null(\mathcal{A})=Span\left\{ \left[ 
\begin{array}{c}
1 \\ 
0 \\ 
{{\mathring{A}}^{-1}(1)} \\ 
0%
\end{array}%
\right] \right\} ,  \label{a_null}
\end{equation}%
where $\mathring{A}:L^{2}(\Omega )\rightarrow L^{2}(\Omega )$ is as given in
(\ref{angst}).

(b) $\sigma _{p}(\mathcal{A}^{\ast }\mathcal{)}\cap \left( i%
\mathbb{R}
\diagdown \{0\}\right) =\varnothing $.
\end{lemma}

\textbf{Proof:} Using the explicit form of the adjoint in Proposition \ref%
{adj}, we can obtain -- in the style of the proof of Proposition \ref{energy}
-- the following relation for all $\Phi =[p_{0},\mathbf{u}%
_{0},w_{1},w_{2}]\in D(\mathcal{A}^{\ast })$:%
\begin{equation}
\func{Re}\left( \mathcal{A}^{\ast }\Phi ,\Phi \right) =-\left( \sigma (%
\mathbf{u}_{0}),\epsilon (\mathbf{u}_{0})\right) _{\mathcal{O}}-\eta
\left\Vert \mathbf{u}_{0}\right\Vert _{\mathcal{O}}^{2}.  \label{a_energy}
\end{equation}%
Subsequently, we can duplicate the proofs of Lemmas \ref{SS} and \ref%
{point_2} to deduce (a) and (b). \ \ \ $\square $

\bigskip

In turn, Lemma \ref{adj_spec}(b) and the classical functional analysis --
see e.g., p. 127 of \cite{friedman} -- yield

\begin{corollary}
\label{residual}Under the assumption that ambient vector field $\mathbf{U}$
satisfies (\ref{extra}), then $\sigma _{r}(\mathcal{A)}\cap i%
\mathbb{R}
=\varnothing $.
\end{corollary}

\medskip

\subsection{The Flow-Structure Semigroup on $[Null(\mathcal{A})]^{\bot }$}

\bigskip

We set 
\begin{equation}
H_{N}=Null(\mathcal{A})\text{; \ \ }H_{N}^{\bot }=[Null(\mathcal{A)}]^{\bot
},  \label{null}
\end{equation}%
where $Null(\mathcal{A})$ and $[Null(\mathcal{A)}]^{\bot }$ are as
characterized respectively in (\ref{N1}) and (\ref{N2}). Then the state
space can be decomposed as 
\begin{equation}
\mathcal{H=}H_{N}\oplus H_{N}^{\bot }.  \label{state}
\end{equation}%
Using either the compatibility condition which characterizes $[Null(\mathcal{%
A)}]^{\bot }$ in (\ref{N2}) -- or alternatively, using the fact that $Null(%
\mathcal{A)=}Null(\mathcal{A}^{\ast }\mathcal{)}$, as expressed in (\ref%
{a_null}) -- one infers the invariances 
\begin{equation*}
\mathcal{A}\mid _{H_{N}^{\bot }}:D(\mathcal{A}\mid _{H_{N}^{\bot
}})\rightarrow H_{N}^{\bot }\,\text{, \ \ }\left. e^{\mathcal{A}%
t}\right\vert _{H_{N}^{\bot }}:H_{N}^{\bot }\rightarrow H_{N}^{\bot }\text{
for all }t\geq 0\text{,}
\end{equation*}%
where $D(\mathcal{A}\mid _{H_{N}^{\bot }})=D(\mathcal{A})\cap H_{N}^{\bot }$%
. Accordingly, concerning the flow-structure PDE (\ref{system1})--(\ref{IC_2}%
), with initial data restricted to $H_{N}^{\bot }$, we have the following
conclusions:\newline

\begin{proposition}
\label{invariant} (Theorem 3.1 of \cite{agw}) Under the assumption that
ambient vector field $\mathbf{U}$ satisfies (\ref{extra}), then $\mathcal{A}%
\mid _{H_{N}^{\bot }}$generates a $C_{0}-$\emph{semigroup of contractions }$%
\exp \left( \left. \mathcal{A}\right\vert _{H_{N}^{\bot }}t\right) =\left.
e^{\mathcal{A}t}\right\vert _{H_{N}^{\bot }}\in \mathcal{L}(H_{N}^{\bot })$.
Consequently,

(a) If initial data $[p_{0},\mathbf{u}_{0},w_{1},w_{2}]\in D(\mathcal{A}\mid
_{H_{N}^{\bot }}),$ then the solution $[p,\mathbf{u},w,w_{t}]$ of (\ref%
{system1})--(\ref{IC_2}) is in \newline
$C([0,\infty );D(\mathcal{A}\mid _{H_{N}^{\bot }}))\cap C^{1}([0,\infty
);H_{N}^{\bot }).$\newline

(b) If initial data $[p_{0},\mathbf{u}_{0},w_{1},w_{2}]\in H_{N}^{\bot },$
then the solution $[p,\mathbf{u},w,w_{t}]$ of (\ref{system1})--(\ref{IC_2})
is in $C([0,\infty );H_{N}^{\bot })$. In addition, there is the dissipative
relation for all $0<t<\infty $, 
\begin{equation*}
\int_{0}^{t}\left[ \left( \sigma (\mathbf{u}(s),\epsilon (\mathbf{u}%
(s))\right) _{\mathcal{O}}+\eta \left\Vert \mathbf{u(s)}\right\Vert _{%
\mathcal{O}}^{2}\right] ds=\left\Vert [p_{0},\mathbf{u}_{0},w_{1},w_{2}]%
\right\Vert _{\mathcal{H}}^{2}-\left\Vert [p(t),\mathbf{u}%
(t),w(t),w_{t}(t)]\right\Vert _{\mathcal{H}}^{2};
\end{equation*}%
whence by Korn's Inequality and the contraction of the semigroup\ $\left\{
\exp \left( \left. \mathcal{A}\right\vert _{H_{N}^{\bot }}t\right) \right\}
_{t\geq 0}$, \newline
we have $\mathbf{u}\in L^{2}((0,\infty );\mathbf{H}^{1}(\mathcal{O})).$ 
\newline
\end{proposition}

\medskip

By the definition of $H_{N}^{\bot }$, the point $\lambda =0$ is not an
eigenvalue of $\mathcal{A}\mid _{H_{N}^{\bot }}.$ In point of fact, $\lambda
=0$ is in the resolvent set of $\mathcal{A}\mid _{H_{N}^{\bot }}$:

\medskip

\begin{lemma}
\label{resolve}The point $\lambda =0$ is in the resolvent set $\rho $ $(%
\mathcal{A}\mid _{H_{N}^{\bot }})$ of $\mathcal{A}\mid _{H_{N}^{\bot }}:D(%
\mathcal{A}\mid _{H_{N}^{\bot }})\rightarrow H_{N}^{\bot }.$ That is, $%
\left( \mathcal{A}\mid _{H_{N}^{\bot }}\right) ^{-1}\in \mathcal{L}%
(H_{N}^{\bot }).$
\end{lemma}

\medskip

\textbf{Proof: }As before, we will use the denotations 
\begin{equation*}
\Phi \equiv 
\begin{bmatrix}
p_{0} \\ 
\mathbf{u}_{0} \\ 
w_{1} \\ 
w_{2}%
\end{bmatrix}%
~~\text{ and}~~\Phi ^{\ast }\equiv 
\begin{bmatrix}
p_{0}^{\ast } \\ 
\mathbf{u}_{0}^{\ast } \\ 
w_{1}^{\ast } \\ 
w_{2}^{\ast }%
\end{bmatrix}%
.
\end{equation*}%
The proof is based upon the use of Lemma \ref{inverse} of the Appendix.
Since $\mathcal{A}\mid _{H_{N}^{\bot }}:D(\mathcal{A}\mid _{H_{N}^{\bot
}})\rightarrow H_{N}^{\bot }$ generates a $C_{0}-$semigroup on $H_{N}^{\bot
},$ then by the Hille-Yosida Theorem, it is closed. Also, since $\lambda =0$
is not an eigenvalue of $\mathcal{A}^{\ast }\mid _{H_{N}^{\bot }}$, by Lemma %
\ref{adj_spec} (a), then $Range(\mathcal{A}\mid _{H_{N}^{\bot }})$ is dense
in $H_{N}^{\bot }.$ Thus, to show that $0\in \rho (\mathcal{A}\mid
_{H_{N}^{\bot }}),$ it will suffice by Lemma \ref{inverse} to establish the
following inequality:%
\begin{equation}
\left\Vert \Phi \right\Vert _{\mathcal{H}}\leq C\left\Vert \mathcal{A}\Phi
\right\Vert _{\mathcal{H}},\text{ \ \ }\forall \text{ }\Phi \in D(\mathcal{A}%
\mid _{H_{N}^{\bot }}).  \label{4.32}
\end{equation}%
Keeping this estimate in mind, we consider the relation%
\begin{equation}
\mathcal{A}\Phi =\Phi ^{\ast },  \label{4.33}
\end{equation}%
for given $\Phi \in D(\mathcal{A}\mid _{H_{N}^{\bot }}).$ In PDE terms, (\ref%
{4.33}) is the following static system:%
\begin{align}
& \left\{ 
\begin{array}{l}
-\mathbf{U}\cdot \nabla p_{0}-div~\mathbf{u}_{0}=p_{0}^{\ast }~\text{ in }~%
\mathcal{O} \\ 
-\nabla p_{0}+div~\sigma (\mathbf{u}_{0})-\eta \mathbf{u}_{0}-\mathbf{U}%
\cdot \nabla \mathbf{u}_{0}=\mathbf{u}_{0}^{\ast }~\text{ in }~\mathcal{O}
\\ 
w_{2}=w_{1}^{\ast }\text{ \ in }~\Omega \\ 
-\Delta ^{2}w_{1}-\left[ 2\nu \partial _{x_{3}}(\mathbf{u}_{0})_{3}+\lambda 
\text{div}(\mathbf{u}_{0})-p_{0}\right] _{\Omega }=w_{2}^{\ast }~\text{ in }%
~\Omega%
\end{array}%
\right.  \label{4.34} \\
&  \notag \\
& \left\{ 
\begin{array}{l}
(\sigma (\mathbf{u}_{0})\mathbf{n}-p_{0}\mathbf{n})\cdot \boldsymbol{\tau }%
=0~\text{ on }~\partial \mathcal{O} \\ 
\mathbf{u}_{0}\cdot \mathbf{n}=0~\text{ on }~S \\ 
\mathbf{u}_{0}\cdot \mathbf{n}=w_{2}~\text{ on }~\Omega \\ 
w_{1}=\frac{\partial w_{1}}{\partial \nu }=0~\text{ on }~\partial \Omega .%
\end{array}%
\right.  \label{4.35}
\end{align}%
Applying Proposition \ref{energy} to (\ref{4.34})-(\ref{4.35}), we have the
dissipative relation,%
\begin{equation}
\left( \sigma (\mathbf{u}_{0}),\epsilon (\mathbf{u}_{0})\right) +\eta
\left\Vert \mathbf{u}_{0}\right\Vert ^{2}=\left\vert \func{Re}\left( \Phi
^{\ast },\Phi \right) _{\mathcal{H}}\right\vert .  \label{4.36}
\end{equation}

\smallskip

In addition, if 
\begin{equation}
p_{0}=q_{0}+c_{0},  \label{4.37}
\end{equation}%
where 
\begin{equation}
\int\limits_{\mathcal{O}}q_{0}d\mathcal{O=}0,~~c_{0}=\text{constant,}
\label{sec}
\end{equation}%
then from estimate (\ref{med2}) of Lemma \ref{med}, we have%
\begin{equation}
\left\Vert q_{0}\right\Vert _{\mathcal{O}}+\left\Vert \sigma (\mathbf{u}_{0})%
\mathbf{n}-p_{0}\mathbf{n}\right\Vert _{\mathbf{H}^{-\frac{1}{2}}(\partial 
\mathcal{O})}\leq C_{\beta }\left( \sqrt{\left\Vert \Phi \right\Vert _{%
\mathcal{H}}\left\Vert \Phi ^{\ast }\right\Vert _{\mathcal{H}}}+\left\Vert
\Phi ^{\ast }\right\Vert _{\mathcal{H}}.\right)  \label{s2}
\end{equation}

\smallskip

Next, if we multiply the mechanical equation in (\ref{4.34}) by $\overline{%
w_{1}}$, integrate, and then integrate by parts, we get -- using the
decomposition (\ref{4.37}), 
\begin{equation}
-\left\Vert \Delta w_{1}\right\Vert _{\Omega }^{2}-\left\langle \sigma (%
\mathbf{u}_{0})\mathbf{n}-q_{0}\mathbf{n,}\left( w_{1}\mathbf{n}\right)
_{ext}\right\rangle _{\partial \mathcal{O}}+\left( c_{0},w_{1}\right)
_{\Omega }=\left( w_{2}^{\ast }~,w_{1}\right) _{\Omega }.  \label{4.40}
\end{equation}%
Here, $\mathbf{H}^{2}-$function $\left( w_{1}\mathbf{n}\right) _{ext}$ is
again given by%
\begin{equation}
\left( w_{1}\mathbf{n}\right) _{ext}(x)=%
\begin{cases}
0 & ~\text{ }x\in \text{ }~S \\ 
w_{1}\mathbf{n} & ~\text{ }x\in ~\Omega .%
\end{cases}
\label{4.41}
\end{equation}%
Now, using the characterization of $H_{N}^{\bot }=[Null(\mathcal{A)}]^{\bot
} $, as given in Lemma \ref{SS}(b), we have that%
\begin{equation}
\left( c_{0},w_{1}\right) =c_{0}\int\limits_{\Omega }\overline{w_{1}}%
=-[meas(\Omega )]\left\vert c_{0}\right\vert ^{2}.  \label{4.42}
\end{equation}%
Applying (\ref{4.42}) to (\ref{4.40}), we then have%
\begin{equation*}
\left\Vert \Delta w_{1}\right\Vert ^{2}+[meas(\Omega
)]c_{0}^{2}=-\left\langle \sigma (\mathbf{u}_{0})\mathbf{n}-q_{0}\mathbf{n,}%
\left( w_{1}\mathbf{n}\right) _{ext}\right\rangle -\left( w_{2}^{\ast
}~,w_{1}\right) ;
\end{equation*}%
Whence we obtain%
\begin{equation*}
(1-\epsilon _{0})\left\Vert \Delta w_{1}\right\Vert ^{2}+[meas(\Omega
)]c_{0}^{2}\leq C_{\epsilon }\left\{ \left\Vert \sigma (\mathbf{u}_{0})%
\mathbf{n}-q_{0}\mathbf{n}\right\Vert _{\mathbf{H}^{-\frac{1}{2}}(\partial 
\mathcal{O})}^{2}+\left\Vert \Phi ^{\ast }\right\Vert _{\mathcal{H}%
}^{2}\right\} .
\end{equation*}%
Applying the estimate (\ref{s2}) to right hand side, as well as Young's
Inequality, we have then%
\begin{equation}
(1-\epsilon _{0})\left\Vert \Delta w_{1}\right\Vert ^{2}+[meas(\Omega
)]c_{0}^{2}\leq \delta \left\Vert \Phi \right\Vert _{\mathcal{H}%
}^{2}+C_{\delta }\left\Vert \Phi ^{\ast }\right\Vert _{\mathcal{H}}^{2}.
\label{4.43}
\end{equation}

\medskip

To conclude the proof: Estimating the normal component boundary condition in
(\ref{4.35}) with the Sobolev Trace Theorem, (\ref{4.36}) and Korn's
Inequality, we have 
\begin{equation}
\left\Vert w_{2}\right\Vert _{\Omega }^{2}\leq \delta \left\Vert \Phi
\right\Vert _{\mathcal{H}}^{2}+C_{\delta }\left\Vert \Phi ^{\ast
}\right\Vert _{\mathcal{H}}^{2}.  \label{4.44}
\end{equation}%
Combining the relations (\ref{4.36}), (\ref{4.37}), (\ref{s2}), (\ref{4.43}%
), and (\ref{4.44}) (and rescaling $\delta >0$), we have finally

\begin{equation*}
\left\Vert p_{0}\right\Vert _{\mathcal{O}}^{2}+\left\Vert \mathbf{u}%
_{0}\right\Vert _{\mathbf{H}^{1}(\mathcal{O})}^{2}+\left\Vert \Delta
w_{1}\right\Vert _{\Omega }^{2}+\left\Vert w_{2}\right\Vert _{\Omega
}^{2}\leq \delta \left\Vert \Phi \right\Vert _{\mathcal{H}}^{2}+C_{\delta
}\left\Vert \Phi ^{\ast }\right\Vert _{\mathcal{H}}^{2}.
\end{equation*}%
Taking $\delta >0$ sufficiently small now gives the estimate (\ref{4.32}),
as required. This concludes the proof of Lemma \ref{resolve}. \ \ \ $\square 
$

\medskip

By combining Lemma \ref{point_2}, Corollary \ref{coro}, Corollary \ref%
{residual}, and Lemma \ref{resolve}, we have now,

\begin{proposition}
\label{no_spec}Concerning $\mathcal{A}\mid _{H_{N}^{\bot }}:D(\mathcal{A}%
\mid _{H_{N}^{\bot }})\rightarrow H_{N}^{\bot }$: Under the assumption that
ambient vector field $\mathbf{U}$ satisfies (\ref{extra}), then $i\mathbb{R}%
\subset \rho (\mathcal{A}\mid _{H_{N}^{\bot }}).$
\end{proposition}

\smallskip

In turn, Propositions \ref{invariant}, \ref{no_spec} and Theorem \ref%
{stab-st} of the Appendix (and in \cite{arendt}) give,

\smallskip

\begin{proposition}
\label{strong}Under the assumption that ambient vector field $\mathbf{U}$
satisfies (\ref{extra}), if initial data $[p_{0},\mathbf{u}%
_{0},w_{1},w_{2}]\in H_{N}^{\bot },$ then the solution $[p,\mathbf{u}%
,w,w_{t}]$ of (\ref{system1})--(\ref{IC_2}) decays asymptotically to the
zero state; i.e., \newline
$\lim_{t\rightarrow \infty }\left\Vert [p(t),\mathbf{u}(t),w(t),w_{t}(t)]%
\right\Vert _{\mathcal{H}}=0$.
\end{proposition}

\medskip

We are now in a position to establish the stronger and main result of
uniform decay of solutions to the compressible flow-structure PDE model (\ref%
{system1})--(\ref{IC_2}).

\section{\label{freq}The Proof of Theorem \protect\ref{uniform}}

The proof of Theorem \ref{uniform} hinges upon an appropriate use of the
stability criterion given in Theorem \ref{stab-un} of the Appendix; see \cite%
{huang} and \cite{pruss}. To this end, and given the fact that the imaginary
axis is contained in the resolvent set of $\mathcal{A}\mid _{H_{N}^{\bot
}}:D(\mathcal{A}\mid _{H_{N}^{\bot }})\rightarrow H_{N}^{\bot }$, by
Proposition \ref{no_spec}, our main goal is to show the following estimate
for the resolvent operator on the imaginary axis:%
\begin{equation}
\left\Vert \left( i\beta I-\mathcal{A}\mid _{H_{N}^{\bot }}\right)
^{-1}\right\Vert _{\mathcal{L}(H_{N}^{\bot })}\leq C^{\ast }\text{, \ for
all }\beta \in \mathbb{R}\text{, and some }C^{\ast }>0,\text{ }  \label{5.1}
\end{equation}%
where the constant $C^{\ast }$ is independent of parameter $\beta \in 
\mathbb{R}
.$ With a view of establishing the estimate (\ref{5.1}), we consider the
equation%
\begin{equation}
\left( i\beta I-\mathcal{A}\right) \Phi =\Phi ^{\ast },  \label{5.2}
\end{equation}%
where 
\begin{equation}
\text{data }\Phi ^{\ast }=%
\begin{bmatrix}
p_{0}^{\ast } \\ 
\mathbf{u}_{0}^{\ast } \\ 
w_{1}^{\ast } \\ 
w_{2}^{\ast }%
\end{bmatrix}%
\in H_{N}^{\bot }\text{, and so corresponding solution }\Phi \equiv 
\begin{bmatrix}
p_{0} \\ 
\mathbf{u}_{0} \\ 
w_{1} \\ 
w_{2}%
\end{bmatrix}%
\in D(\mathcal{A}\mid _{H_{N}^{\bot }}).  \label{image}
\end{equation}

\bigskip

In PDE terms, the resolvent relation (\ref{5.2}) becomes%
\begin{align}
& \left\{ 
\begin{array}{l}
i\beta p_{0}+\mathbf{U}\cdot \nabla p_{0}+\func{div}\mathbf{u}%
_{0}=p_{0}^{\ast }~\text{ in }~\mathcal{O} \\ 
i\beta \mathbf{u}_{0}+\nabla p_{0}-\func{div}\sigma (\mathbf{u}_{0})+\eta 
\mathbf{u}_{0}+\mathbf{U}\cdot \nabla \mathbf{u}_{0}=\mathbf{u}_{0}^{\ast }~%
\text{ in }~\mathcal{O} \\ 
i\beta w_{1}-w_{2}=w_{1}^{\ast }\text{ \ in }~\Omega \\ 
i\beta w_{2}+\Delta ^{2}w_{1}+\left[ 2\nu \partial _{x_{3}}(\mathbf{u}%
_{0})_{3}+\lambda \text{div}(\mathbf{u}_{0})-p_{0}\right] _{\Omega
}=w_{2}^{\ast }~\text{ in }~\Omega%
\end{array}%
\right.  \label{5.3} \\
&  \notag \\
& \left\{ 
\begin{array}{l}
(\sigma (\mathbf{u}_{0})\mathbf{n}-p_{0}\mathbf{n})\cdot \boldsymbol{\tau }%
=0~\text{ on }~\partial \mathcal{O} \\ 
\mathbf{u}_{0}\cdot \mathbf{n}=0~\text{ on }~S \\ 
\mathbf{u}_{0}\cdot \mathbf{n}=w_{2}~\text{ on }~\Omega \\ 
w_{1}=\frac{\partial w_{1}}{\partial \nu }=0~\text{ on }~\partial \Omega .%
\end{array}%
\right.  \label{5.4}
\end{align}

By way of proving the estimate (\ref{5.1}), \ we start with the following
lemma:

\begin{lemma}
\label{US-L} With $L^{2}$-pressure component $p_{0}$ of solution $\Phi $ of (%
\ref{5.2}) enjoying the decomposition 
\begin{equation}
p_{0}=q_{0}+c_{0},  \label{pre}
\end{equation}%
where 
\begin{equation}
L^{2}\text{-function }q_{0}\text{ satisfies }\int\limits_{\mathcal{O}}q_{0}d%
\mathcal{O=}0\text{, and }c_{0}=\text{constant;}  \label{d_2}
\end{equation}%
then the following estimate holds for $\Phi $: 
\begin{equation}
\begin{array}{l}
\left\Vert \mathbf{u}_{0}\right\Vert _{\mathbf{H}^{1}(\mathcal{O}%
)}^{2}+\left\Vert w_{2}\right\Vert _{\Omega }^{2}\leq C\left\vert \left(
\Phi ^{\ast },\Phi \right) _{\mathcal{H}}\right\vert ; \\ 
\\ 
meas(\Omega )\left\vert c_{0}\right\vert ^{2}+\frac{1}{2}\left\Vert \Delta
w_{1}\right\Vert _{\Omega }^{2}\leq C\left( \left\vert \left( \Phi ^{\ast
},\Phi \right) _{\mathcal{H}}\right\vert +\left\Vert \Phi ^{\ast
}\right\Vert _{\mathcal{H}}^{2}\right) .%
\end{array}
\label{start}
\end{equation}
\end{lemma}

\vspace{0.5cm} \textbf{Proof of Lemma \ref{US-L}}. Initially, we proceed as
before: Taking the $\mathcal{H}$-inner product of both sides of (\ref{5.2})
with respect to $\Phi $, and subsequently appealing to Proposition \ref%
{energy}, we have then%
\begin{equation*}
\left( \sigma (\mathbf{u}),\epsilon (\mathbf{u})\right) _{\mathcal{O}}+\eta
\left\Vert \mathbf{u}\right\Vert _{\mathcal{O}}^{2}=\func{Re}\left( \Phi
^{\ast },\Phi \right) _{\mathcal{H}},
\end{equation*}%
whence by Korn's Inequality we have%
\begin{equation}
\left\Vert \mathbf{u}_{0}\right\Vert _{\mathbf{H}^{1}(\mathcal{O})}^{2}\leq
C\left\vert \left( \Phi ^{\ast },\Phi \right) _{\mathcal{H}}\right\vert .
\label{5.5}
\end{equation}%
In turn, this estimate, the normal component boundary condition in (\ref{5.4}%
), and the Sobolev Trace Theorem give%
\begin{equation}
\left\Vert w_{2}\right\Vert _{\Omega }^{2}\leq C\left\vert \left( \Phi
^{\ast },\Phi \right) _{\mathcal{H}}\right\vert .  \label{5.55}
\end{equation}
Moreover, appealing to Lemma \ref{med}, with pressure component $q_{0}$ as
given in (\ref{d_2}), we have the estimate%
\begin{equation}
\left\Vert \sigma (\mathbf{u}_{0})\mathbf{n}-q_{0}\mathbf{n}\right\Vert _{%
\mathbf{H}^{-\frac{1}{2}}(\partial \mathcal{O})}\leq C\left( \sqrt{%
\left\Vert \Phi \right\Vert _{\mathcal{H}}\left\Vert \Phi ^{\ast
}\right\Vert _{\mathcal{H}}}+\left\Vert \Phi ^{\ast }\right\Vert _{\mathcal{H%
}}+\left\vert \beta \right\vert \left\Vert \mathbf{u}_{0}\right\Vert _{%
\mathcal{O}}\right) .  \label{5.6}
\end{equation}

Next, we multiply the structural equation in (\ref{5.3}) by $\overline{w_{1}}
$, integrate, and then integrate by parts -- and use the decomposition in (%
\ref{d_2}) -- to have%
\begin{equation}
\left\Vert \Delta w_{1}\right\Vert _{\Omega }^{2}=-\left\langle \sigma (%
\mathbf{u}_{0})\mathbf{n}-q_{0}\mathbf{n,}\left( w_{1}\mathbf{n}\right)
_{ext}\right\rangle _{\partial \mathcal{O}}+c_{0}\int\limits_{\Omega }%
\overline{w_{1}}d\Omega -i\beta \left( w_{2},w_{1}\right) _{\Omega }+\left(
w_{2}^{\ast }~,w_{1}\right) _{\Omega };  \label{5.10}
\end{equation}%
where as before $\mathbf{H}^{2}-$function $\left( w_{1}\mathbf{n}\right)
_{ext}$ is given by%
\begin{equation*}
\left( w_{1}\mathbf{n}\right) _{ext}(x)=%
\begin{cases}
0 & ~\text{ }x\in \text{ }~S \\ 
w_{1}\mathbf{n} & ~\text{ }x\in ~\Omega .%
\end{cases}%
\end{equation*}%
Applying now the compatibility condition for $H_{N}^{\bot }=[Null(\mathcal{A}%
)]^{\bot }$ in Lemma \ref{SS}(b), we thus obtain%
\begin{equation}
meas(\Omega )\left\vert c_{0}\right\vert ^{2}+\left\Vert \Delta
w_{1}\right\Vert _{\Omega }^{2}=-\left\langle \sigma (\mathbf{u}_{0})\mathbf{%
n}-q_{0}\mathbf{n,}\left( w_{1}\mathbf{n}\right) _{ext}\right\rangle
_{\partial \mathcal{O}}-i\beta \left( w_{2},w_{1}\right) _{\Omega }+\left(
w_{2}^{\ast }~,w_{1}\right) _{\Omega }.  \label{5.13}
\end{equation}

\smallskip

We now consider the respective cases $\left\vert \beta \right\vert >1$ and $%
\left\vert \beta \right\vert \leq 1$.\newline

\medskip

\underline{\textbf{Case I}} \textbf{(}$\left\vert \beta \right\vert >1$%
\textbf{:}\newline
Using the third resolvent relation in (\ref{5.3}), we have that 
\begin{equation*}
w_{1}=-\frac{i}{\beta }w_{2}-\frac{i}{\beta }w_{1}^{\ast }.
\end{equation*}%
Accordingly the first term on RHS of (\ref{5.13}) becomes%
\begin{equation}
-\left\langle \sigma (\mathbf{u}_{0})\mathbf{n}-q_{0}\mathbf{n,}\left( w_{1}%
\mathbf{n}\right) _{ext}\right\rangle _{\partial \mathcal{O}}=-\frac{i}{%
\beta }\left\langle \sigma (\mathbf{u}_{0})\mathbf{n}-q_{0}\mathbf{n,}\left(
w_{2}\mathbf{n}\right) _{ext}\right\rangle _{\partial \mathcal{O}}-\frac{i}{%
\beta }\left\langle \sigma (\mathbf{u}_{0})\mathbf{n}-q_{0}\mathbf{n,}\left(
w_{1}^{\ast }\mathbf{n}\right) _{ext}\right\rangle _{\partial \mathcal{O}}.
\label{5.14}
\end{equation}

\smallskip

\textbf{(a)} To handle the first term on RHS of (\ref{5.14}): we use the
tangential B.C. in (\ref{5.4}), and the fact that $\Phi \in D(\mathcal{A}%
\mid _{H_{N}^{\bot }})\,$, which means that flow solution component $\mathbf{%
u}_{0}$ on the boundary has the decomposition%
\begin{equation}
\mathbf{u}_{0}\mid _{\partial \mathcal{O}}=\mathbf{f}_{0}\mid _{\partial 
\mathcal{O}}+\left( w_{2}\mathbf{n}\right) _{ext},  \label{5.145}
\end{equation}%
where $\mathbf{f}_{0}\in \mathbf{V}_{0}$ \noindent (and so $\mathbf{f}%
_{0}\mid _{\partial \mathcal{O}}\in TH^{1/2}(\partial \mathcal{O})$).

\medskip

In addition, given that pressure component $p_{0}$ of solution $\Phi \in D(%
\mathcal{A}\mid _{H_{N}^{\bot }})$ satisfies the relation 
\begin{equation}
\left\langle \sigma (\mathbf{u}_{0})\mathbf{n}-p_{0}\mathbf{n,}\tau
\right\rangle _{\partial \mathcal{O}}=0\text{ \ for every }\tau \in
TH^{1/2}(\partial \mathcal{O}),  \label{5.146}
\end{equation}%
by domain requirement (A.iv), then so too does $L^{2}$-component $q_{0}$ of (%
\ref{pre}). \ Indeed, given $\tau \in TH^{1/2}(\partial \mathcal{O})$, let $%
\tilde{\mu}\in \mathbf{H}^{1}(\mathcal{O})$ satisfy the boundary value
problem%
\begin{eqnarray}
\func{div}(\tilde{\mu}) &=&0\text{ \ in }\mathcal{O}  \notag \\
\left. \tilde{\mu}\right\vert _{\partial \mathcal{O}} &=&\tau \text{ \ on }%
\partial \mathcal{O}.  \label{bvp}
\end{eqnarray}%
The existence of such a function $\tilde{\mu}(\tau )$ is assured; see e.g.,
p. 127 of \cite{galdi}. Therewith, we have by (\ref{pre}) and Green's formula%
\begin{eqnarray}
\left\langle \sigma (\mathbf{u}_{0})\mathbf{n}-q_{0}\mathbf{n,}\tau
\right\rangle _{\partial \mathcal{O}} &=&c_{0}\left\langle \mathbf{n,}\tau
\right\rangle _{\partial \mathcal{O}}+\left\langle \sigma (\mathbf{u}_{0})%
\mathbf{n}-p_{0}\mathbf{n,}\tau \right\rangle _{\partial \mathcal{O}}  \notag
\\
&=&(c_{0},\func{div}(\tilde{\mu}(\tau )))_{\mathcal{O}}+(\nabla c_{0},\tilde{%
\mu}(\tau ))_{\mathcal{O}}+0  \notag \\
&=&0\text{, \ for every }\tau \in TH^{1/2}(\partial \mathcal{O}).  \label{Gf}
\end{eqnarray}

Using this relation and (\ref{5.145}) we now obtain%
\begin{eqnarray*}
\left\langle \sigma (\mathbf{u}_{0})\mathbf{n}-q_{0}\mathbf{n,}\left( w_{2}%
\mathbf{n}\right) _{ext}\right\rangle _{\partial \mathcal{O}}
&=&\left\langle \sigma (\mathbf{u}_{0})\mathbf{n}-q_{0}\mathbf{n,f}%
_{0}+\left( w_{2}\mathbf{n}\right) _{ext}\right\rangle _{\partial \mathcal{O}%
} \\
&=&\left\langle \sigma (\mathbf{u}_{0})\mathbf{n}-q_{0}\mathbf{n,u}_{0}\mid
_{\partial \mathcal{O}}\right\rangle ,
\end{eqnarray*}%
and so%
\begin{equation*}
\left\vert \left\langle \sigma (\mathbf{u}_{0})\mathbf{n}-q_{0}\mathbf{n,}%
\left( w_{2}\mathbf{n}\right) _{ext}\right\rangle _{\partial \mathcal{O}%
}\right\vert \leq \left\Vert \sigma (\mathbf{u})\mathbf{n}-q_{0}\mathbf{n}%
\right\Vert _{\mathbf{H}^{-\frac{1}{2}}(\partial \mathcal{O})}\left\Vert 
\mathbf{u}_{0}\mid _{\partial \mathcal{O}}\right\Vert _{H^{\frac{1}{2}%
}(\partial \mathcal{O})}
\end{equation*}%
\begin{equation}
\leq C\left( \left\vert \beta \right\vert \left\Vert \mathbf{u}%
_{0}\right\Vert _{\mathcal{O}}+\sqrt{\left\vert \left( \Phi ^{\ast },\Phi
\right) _{\mathcal{H}}\right\vert }+\left\Vert \Phi ^{\ast }\right\Vert _{%
\mathcal{H}}\right) \sqrt{\left\vert \left( \Phi ^{\ast },\Phi \right) _{%
\mathcal{H}}\right\vert },  \label{5.15}
\end{equation}%
after using (\ref{5.5}) and (\ref{5.6}). Applying once more (\ref{5.5}) to
this estimate, we then have for $\left\vert \beta \right\vert >1$,%
\begin{equation}
\left\vert \frac{i}{\beta }\left\langle \sigma (\mathbf{u})\mathbf{n}-q_{0}%
\mathbf{n,}\left( w_{2}\mathbf{n}\right) _{ext}\right\rangle _{\partial 
\mathcal{O}}\right\vert \leq C\left( \left\vert \left( \Phi ^{\ast },\Phi
\right) _{\mathcal{H}}\right\vert +\left\Vert \Phi ^{\ast }\right\Vert _{%
\mathcal{H}}^{2}\right) ,  \label{5.16}
\end{equation}%
where positive constant $C$ is independent of $\beta $.

\medskip

\textbf{(b) }Now, for the second term on the RHS of (\ref{5.14}): using the
estimate (\ref{5.6}) and (\ref{5.5}) in sequence, we have for $\left\vert
\beta \right\vert >1$, 
\begin{equation*}
\left\vert \frac{i}{\beta }\left\langle \sigma (\mathbf{u})\mathbf{n}-q_{0}%
\mathbf{n,}\left( w_{1}^{\ast }\mathbf{n}\right) _{ext}\right\rangle
_{\partial \mathcal{O}}\right\vert \leq \frac{1}{\left\vert \beta
\right\vert }\left\Vert \sigma (\mathbf{u})\mathbf{n}-q_{0}\mathbf{n}%
\right\Vert _{\mathbf{H}^{-\frac{1}{2}}(\partial \mathcal{O})}\left\Vert
\Delta w_{1}^{\ast }\right\Vert _{\Omega }
\end{equation*}%
\begin{equation*}
\leq C\left( \left\Vert \mathbf{u}_{0}\right\Vert +\sqrt{\left\vert \left(
\Phi ^{\ast },\Phi \right) _{\mathcal{H}}\right\vert }+\left\Vert \Phi
^{\ast }\right\Vert _{\mathcal{H}}\right) \left\Vert \Delta w_{1}^{\ast
}\right\Vert _{\Omega }
\end{equation*}%
\begin{equation}
\leq C\left( \left\vert \left( \Phi ^{\ast },\Phi \right) _{\mathcal{H}%
}\right\vert +\left\Vert \Phi ^{\ast }\right\Vert _{\mathcal{H}}^{2}\right) .
\label{5.17}
\end{equation}

\medskip

Applying now estimates (\ref{5.16}) and (\ref{5.17}) to the RHS of (\ref%
{5.14}), we get%
\begin{equation}
\left\vert \left\langle \sigma (\mathbf{u})\mathbf{n}-q_{0}\mathbf{n,}\left(
w_{1}\mathbf{n}\right) _{ext}\right\rangle _{\partial \mathcal{O}%
}\right\vert \leq C[\left\vert \left( \Phi ^{\ast },\Phi \right) _{\mathcal{H%
}}\right\vert +\left\Vert \Phi ^{\ast }\right\Vert _{\mathcal{H}}^{2}].
\label{5.18}
\end{equation}

\bigskip

We proceed now to the second term on RHS of (\ref{5.13}). For this, we again
use the third resolvent relation in (\ref{5.3}) to have%
\begin{equation*}
\left\vert \beta \left( w_{2},w_{1}\right) _{\Omega }\right\vert =\left\vert
\left( w_{2},w_{2}+w_{1}^{\ast }\right) _{\Omega }\right\vert
\end{equation*}%
\begin{equation*}
\leq \left\Vert w_{2}\right\Vert _{\Omega }^{2}+\left\vert \left(
w_{2},w_{1}^{\ast }\right) _{\Omega }\right\vert
\end{equation*}%
\begin{equation}
\leq C[\left\vert \left( \Phi ^{\ast },\Phi \right) _{\mathcal{H}%
}\right\vert +\left\Vert \Phi ^{\ast }\right\Vert _{\mathcal{H}}^{2}],
\label{5.19}
\end{equation}%
after invoking (\ref{5.55}).

\smallskip

For the third term on RHS of (\ref{5.13}): in similar fashion, we use the
third resolvent relation in (\ref{5.3}), and then (\ref{5.55}), so as to
obtain for $\left\vert \beta \right\vert >1$ 
\begin{equation*}
\left\vert \left( w_{2}^{\ast },w_{1}\right) \right\vert =\left\vert \left(
w_{2}^{\ast },i\beta ^{-1}(w_{2}+w_{1}^{\ast }\right) )\right\vert
\end{equation*}%
\begin{equation}
\leq C[\left\vert \left( \Phi ^{\ast },\Phi \right) _{\mathcal{H}%
}\right\vert +\left\Vert \Phi ^{\ast }\right\Vert _{\mathcal{H}}^{2}].
\label{5.20}
\end{equation}

\smallskip

Applying finally the estimates (\ref{5.18}), (\ref{5.19}), and (\ref{5.20})
to RHS of (\ref{5.13}), we obtain the following estimate for the case $%
\left\vert \beta \right\vert >1$, 
\begin{equation}
meas(\Omega )\left\vert c_{0}\right\vert ^{2}+\left\Vert \Delta
w_{1}\right\Vert _{\Omega }^{2}\leq C\left( \left\vert \left( \Phi ^{\ast
},\Phi \right) _{\mathcal{H}}\right\vert +\left\Vert \Phi ^{\ast
}\right\Vert _{\mathcal{H}}^{2}\right) ,  \label{5.21}
\end{equation}

\bigskip

\underline{\textbf{Case II}} \textbf{(}$\left\vert \beta \right\vert \leq 1$%
\textbf{):}\newline
To estimate the RHS of (\ref{5.13}) in this case, we use estimates (\ref{5.5}%
)-(\ref{5.6}) to have for $\left\vert \beta \right\vert \leq 1$,%
\begin{equation*}
meas(\Omega )\left\vert c_{0}\right\vert ^{2}+\left\Vert \Delta
w_{1}\right\Vert ^{2}=\left\vert \left\langle \sigma (\mathbf{u})\mathbf{n}%
-q_{0}\mathbf{n,}\left( w_{1}\mathbf{n}\right) _{ext}\right\rangle
\right\vert +\left\vert \beta \right\vert \left\vert \left(
w_{2},w_{1}\right) \right\vert +\left\vert \left( w_{2}^{\ast
}~,w_{1}\right) \right\vert
\end{equation*}%
\begin{equation*}
\leq \left\Vert \sigma (\mathbf{u})\mathbf{n}-q_{0}\mathbf{n}\right\Vert _{%
\mathbf{H}^{-\frac{1}{2}}(\partial \mathcal{O})}\left\Vert \Delta
w_{1}\right\Vert _{\Omega }+\left\Vert w_{2}\right\Vert _{\Omega }\left\Vert
w_{1}\right\Vert _{\Omega }+\left\Vert w_{2}^{\ast }\right\Vert _{\Omega
}\left\Vert w_{1}\right\Vert _{\Omega }
\end{equation*}%
\begin{eqnarray*}
&\leq &C\left( \sqrt{\left\Vert \Phi \right\Vert _{\mathcal{H}}\left\Vert
\Phi ^{\ast }\right\Vert _{\mathcal{H}}}+\left\Vert \Phi ^{\ast }\right\Vert
_{\mathcal{H}}+\left\vert \beta \right\vert \left\Vert \mathbf{u}%
_{0}\right\Vert _{\mathcal{O}}\right) .\left\Vert \Delta w_{1}\right\Vert
_{\Omega }+\left\Vert w_{2}^{\ast }\right\Vert _{\Omega }\left\Vert
w_{1}\right\Vert _{\Omega } \\
&\leq &C\left( \sqrt{\left\Vert \Phi \right\Vert _{\mathcal{H}}\left\Vert
\Phi ^{\ast }\right\Vert _{\mathcal{H}}}+\left\Vert \Phi ^{\ast }\right\Vert
_{\mathcal{H}}\right) .\left\Vert \Delta w_{1}\right\Vert _{\Omega }.
\end{eqnarray*}%
An application of Young's inequality with $\delta <\frac{1}{2}$ now gives
that in the case $\left\vert \beta \right\vert \leq 1$%
\begin{equation}
meas(\Omega )\left\vert c_{0}\right\vert ^{2}+\left( 1-\delta \right)
\left\Vert \Delta w_{1}\right\Vert _{\Omega }^{2}\leq C_{\delta }\left[
\left\vert \left( \Phi ^{\ast },\Phi \right) _{\mathcal{H}}\right\vert
+\left\Vert \Phi ^{\ast }\right\Vert _{\mathcal{H}}^{2}\right] .
\label{5.22}
\end{equation}

\smallskip

This finishes the proof of Lemma \ref{US-L}. \ \ \ $\square $

\bigskip

To complete the proof of Theorem \ref{uniform}, we must obtain a beneficial
estimate for the the \textquotedblleft zero average value\textquotedblright\
component $q_{0}$ of pressure term $p_{0}$, as given in (\ref{pre}). For as
things stand now, we only have, by Lemma \ref{med}, the estimate 
\begin{equation}
\left\Vert q_{0}\right\Vert _{\mathcal{O}}\leq C\left( \sqrt{\left\Vert \Phi
\right\Vert _{\mathcal{H}}\left\Vert \Phi ^{\ast }\right\Vert _{\mathcal{H}}}%
+\left\Vert \Phi ^{\ast }\right\Vert _{\mathcal{H}}+\left\vert \beta
\right\vert \left\Vert \mathbf{u}_{0}\right\Vert _{\mathcal{O}}\right) ,
\label{insuf}
\end{equation}

which is unhelpful for $\left\vert \beta \right\vert $ large.

\medskip

By way of refining our estimation of $p_{0}$, we consider the two cases,
\textquotedblleft $\beta $ large\textquotedblright\ and \textquotedblleft $%
\beta $ small\textquotedblright\ in the resolvent equation (\ref{5.2}) (one
would guess that the former is the problematic case).

\bigskip

\textbf{Case I (}$\left\vert \beta \right\vert >1$\textbf{):} We consider
the decomposition%
\begin{equation}
\left\{ 
\begin{array}{l}
q_{0}=\mathfrak{q}_{1}+\mathfrak{q}_{2} \\ 
\mathbf{u}_{0}=\boldsymbol{u}_{1}+\boldsymbol{u}_{2}%
\end{array}%
\right.  \label{5.23}
\end{equation}%
where $(\boldsymbol{u}_{1},\mathfrak{q}_{1})$ solves 
\begin{equation}
\left\{ 
\begin{array}{l}
\nabla \mathfrak{q}_{1}-\func{div}\sigma (\boldsymbol{u}_{1})+\eta 
\boldsymbol{u}_{1}=-i\beta \mathbf{u}_{0}~~~\text{ in }~~~\mathcal{O} \\ 
\begin{array}{l}
~\func{div}(\boldsymbol{u}_{1})=0\text{ in }\Omega ~ \\ 
\boldsymbol{u}_{1}=0\text{ \ on }\partial \mathcal{O,}%
\end{array}%
\end{array}%
\right.  \label{5.24}
\end{equation}%
and $(\boldsymbol{u}_{2},\mathfrak{q}_{2})$ solves%
\begin{equation}
\left\{ 
\begin{array}{l}
\nabla \mathfrak{q}_{2}-\func{div}\sigma (\boldsymbol{u}_{2})+\eta 
\boldsymbol{u}_{2}=\mathbf{-U}\cdot \nabla \mathbf{u}_{0}+\mathbf{u}%
_{0}^{\ast }~~~\text{ in }~~~\mathcal{O} \\ 
\begin{array}{l}
~\func{div}(\boldsymbol{u}_{2})=~\func{div}(\mathbf{u}_{0})\text{ \ in }%
\mathcal{O}~ \\ 
\boldsymbol{u}_{2}=\mathbf{u}_{0}\text{on \ }\partial \mathcal{O}.%
\end{array}%
\end{array}%
\right.  \label{5.25}
\end{equation}%
We apply the results in \cite{dauge} to the PDE component (\ref{5.24}) --
see in particular (1.4)-(1.5), (1.8), (1.9) of \cite{dauge}) with $s=1.%
\mathcal{\ }$This gives 
\begin{equation}
\left\Vert \boldsymbol{u}_{1}\right\Vert _{\mathbf{H}^{2}(\mathcal{O}%
)}+\left\Vert \mathfrak{q}_{1}\right\Vert _{H^{1}(\mathcal{O})}\leq
C\left\vert \beta \right\vert \left\Vert \mathbf{u}_{0}\right\Vert _{%
\mathcal{O}}.  \label{5.26}
\end{equation}

\smallskip

Moreover, applying Theorem \ref{static} of the Appendix to PDE (\ref{5.25}),
followed by Lemma \ref{US-L} and the Sobolev Trace Theorem, we obtain the
estimate%
\begin{eqnarray}
\left\Vert \boldsymbol{u}_{2}\right\Vert _{\mathbf{H}^{1}(\mathcal{O}%
)}+\left\Vert \mathfrak{q}_{2}\right\Vert _{\mathcal{O}} &\leq &C\left[
\left\Vert \mathbf{U\cdot }\nabla \mathbf{u}_{0}+\mathbf{u}_{0}^{\ast
}\right\Vert _{\mathcal{O}}+\left\Vert \func{div}(\mathbf{u}_{0})\right\Vert
_{\mathcal{O}}+\left\Vert \mathbf{u}_{0}\mid _{\partial \mathcal{O}%
}\right\Vert _{H^{\frac{1}{2}}(\partial \mathcal{O})}\right]  \notag \\
&\leq &C\left( \sqrt{\left\vert \left( \Phi ^{\ast },\Phi \right) _{\mathcal{%
H}}\right\vert }+\left\Vert \Phi ^{\ast }\right\Vert _{\mathcal{H}}\right) .
\label{5.266}
\end{eqnarray}

\medskip

With these two estimates in hand, we now consider the following term:%
\begin{equation}
\left( \mathbf{U}\cdot \nabla p_{0},p_{0}\right) _{\mathcal{O}}=\left( 
\mathbf{U}\cdot \nabla p_{0},q_{0}\right) _{\mathcal{O}}+\left( \mathbf{U}%
\cdot \nabla p_{0},c_{0}\right) _{\mathcal{O}},  \label{5.267}
\end{equation}%
where $p_{0}=q_{0}+c_{0}$ is the decomposition (\ref{pre}). Since $\left( 
\mathbf{U}\cdot \nabla p_{0},c_{0}\right) _{\mathcal{O}}=0$ by Green's
Theorem (as $\mathbf{U}$ is divergence free), then from (\ref{5.23}) we have%
\begin{equation}
\left( \mathbf{U}\cdot \nabla p_{0},p_{0}\right) _{\mathcal{O}}=\left( 
\mathbf{U}\cdot \nabla p_{0},\mathfrak{q}_{1}\right) _{\mathcal{O}}+\left( 
\mathbf{U}\cdot \nabla p_{0},\mathfrak{q}_{2}\right) _{\mathcal{O}}.
\label{5.28}
\end{equation}

\medskip

For the first term on RHS of (\ref{5.28}), we have from Green's Theorem and (%
\ref{5.26}), 
\begin{equation*}
\left( \mathbf{U}\cdot \nabla p_{0},\mathfrak{q}_{1}\right) _{\mathcal{O}%
}=-\left( p_{0},\mathbf{U}\cdot \nabla \mathfrak{q}_{1}\right) _{\mathcal{O}%
}\leq C\left\vert \beta \right\vert \left\Vert \mathbf{u}_{0}\right\Vert _{%
\mathcal{O}}\left\Vert p_{0}\right\Vert _{\mathcal{O}}.
\end{equation*}%
Refining this RHS by means of Lemma \ref{US-L}, we then have%
\begin{equation}
\left\vert \left( \mathbf{U}\cdot \nabla p_{0},\mathfrak{q}_{1}\right) _{%
\mathcal{O}}\right\vert \leq C\left\vert \beta \right\vert \left\Vert
p_{0}\right\Vert _{\mathcal{O}}\sqrt{\left\vert \left( \Phi ^{\ast },\Phi
\right) _{\mathcal{H}}\right\vert }.  \label{5.29}
\end{equation}

\medskip

For the second term on RHS of (\ref{5.28}): We write

\begin{equation}
\left( \mathbf{U}\cdot \nabla p_{0},\mathfrak{q}_{2}\right) _{\mathcal{O}%
}=\left( \mathbf{U}\cdot \nabla p_{0}+\func{div}(\mathbf{u}_{0}),\mathfrak{q}%
_{2}\right) _{\mathcal{O}}-\left( \func{div}(\mathbf{u}_{0}),\mathfrak{q}%
_{2}\right) _{\mathcal{O}}.  \label{5.30}
\end{equation}

\medskip

To estimate the first term on RHS of (\ref{5.30}), we use the pressure
equation in (\ref{5.3}), so as to get%
\begin{eqnarray*}
\left\vert \left( \mathbf{U}\cdot \nabla p_{0}+\func{div}(\mathbf{u}_{0}),%
\mathfrak{q}_{2}\right) _{\mathcal{O}}\right\vert &\leq &\left\Vert \mathbf{U%
}\cdot \nabla p_{0}+\func{div}(\mathbf{u}_{0})\right\Vert _{\mathcal{O}%
}\left\Vert \mathfrak{q}_{2}\right\Vert _{\mathcal{O}} \\
&=&\left\Vert i\beta p_{0}-p_{0}^{\ast }\right\Vert _{\mathcal{O}}\left\Vert 
\mathfrak{q}_{2}\right\Vert _{\mathcal{O}}.
\end{eqnarray*}%
We have then, after using (\ref{5.266}) and Young's Inequality,%
\begin{eqnarray}
\left\vert \left( \mathbf{U}\cdot \nabla p_{0}+\func{div}(\mathbf{u}_{0}),%
\mathfrak{q}_{2}\right) _{\mathcal{O}}\right\vert &\leq &C\left\vert \beta
\right\vert \left\Vert p_{0}\right\Vert _{\mathcal{O}}\left( \sqrt{%
\left\vert \left( \Phi ^{\ast },\Phi \right) _{\mathcal{H}}\right\vert }%
+\left\Vert \Phi ^{\ast }\right\Vert _{\mathcal{H}}\right)  \notag \\
&&+C\left\Vert p_{0}^{\ast }\right\Vert _{\mathcal{O}}\left( \sqrt{%
\left\vert \left( \Phi ^{\ast },\Phi \right) _{\mathcal{H}}\right\vert }%
+\left\Vert \Phi ^{\ast }\right\Vert _{\mathcal{H}}\right)  \notag \\
&&  \notag \\
&\leq &C\left\vert \beta \right\vert \left\Vert p_{0}\right\Vert _{\mathcal{O%
}}\left( \sqrt{\left\vert \left( \Phi ^{\ast },\Phi \right) _{\mathcal{H}%
}\right\vert }+\left\Vert \Phi ^{\ast }\right\Vert _{\mathcal{H}}\right)
+\epsilon \left\Vert \Phi \right\Vert _{\mathcal{H}}^{2}+C_{\epsilon
}\left\Vert \Phi ^{\ast }\right\Vert _{\mathcal{H}}^{2}.  \label{5.31}
\end{eqnarray}

\medskip

For the second term on RHS of (\ref{5.30}), we use Lemma \ref{US-L} and (\ref%
{5.266}) (and Young's Inequality) to have%
\begin{eqnarray*}
\left\vert \left( \func{div}(\mathbf{u}_{0}),\mathfrak{q}_{2}\right)
\right\vert &\leq &C\sqrt{\left\vert \left( \Phi ^{\ast },\Phi \right) _{%
\mathcal{H}}\right\vert }\left( \sqrt{\left\vert \left( \Phi ^{\ast },\Phi
\right) _{\mathcal{H}}\right\vert }+\left\Vert \Phi ^{\ast }\right\Vert _{%
\mathcal{H}}\right) \\
&\leq &\epsilon \left\Vert \Phi \right\Vert _{\mathcal{H}}^{2}+C_{\epsilon
}\left\Vert \Phi ^{\ast }\right\Vert _{\mathcal{H}}^{2}\text{.}
\end{eqnarray*}

This estimate gives then, in combination with (\ref{5.30}) and (\ref{5.31})
(and a rescaling of $\epsilon >0$),%
\begin{equation}
\left( \mathbf{U}\cdot \nabla p_{0},\mathfrak{q}_{2}\right) \leq C\left\vert
\beta \right\vert \left\Vert p_{0}\right\Vert _{\mathcal{O}}\left( \sqrt{%
\left\vert \left( \Phi ^{\ast },\Phi \right) _{\mathcal{H}}\right\vert }%
+\left\Vert \Phi ^{\ast }\right\Vert _{\mathcal{H}}\right) +\epsilon
\left\Vert \Phi \right\Vert _{\mathcal{H}}^{2}+C_{\epsilon }\left\Vert \Phi
^{\ast }\right\Vert _{\mathcal{H}}^{2}.  \label{5.33}
\end{equation}

\medskip

Combining (\ref{5.28}), (\ref{5.29}) and (\ref{5.33}), we now have

\begin{equation}
\left\vert \left( \mathbf{U}\cdot \nabla p_{0},p_{0}\right) \right\vert \leq
C\left\vert \beta \right\vert \left\Vert p_{0}\right\Vert _{\mathcal{O}%
}\left( \sqrt{\left\vert \left( \Phi ^{\ast },\Phi \right) _{\mathcal{H}%
}\right\vert }+\left\Vert \Phi ^{\ast }\right\Vert _{\mathcal{H}}\right)
+\epsilon \left\Vert \Phi \right\Vert _{\mathcal{H}}^{2}+C_{\epsilon
}\left\Vert \Phi ^{\ast }\right\Vert _{\mathcal{H}}^{2}.  \label{5.34}
\end{equation}

\medskip

In order to use this estimate, we multiply both sides of the pressure
equation in (\ref{5.3}) by $\overline{p_{0}}$, and integrate: This gives%
\begin{equation}
i\beta \left\Vert p_{0}\right\Vert _{\mathcal{O}}^{2}=-\left( \mathbf{U}%
\cdot \nabla p_{0},p_{0}\right) _{\mathcal{O}}-\left( \func{div}(\mathbf{u}%
_{0}),p_{0}\right) _{\mathcal{O}}+\left( {p}_{0}^{\ast },p_{0}\right) _{%
\mathcal{O}}.  \label{5.35}
\end{equation}%
Estimating the RHS by (\ref{5.34}) and Lemma \ref{US-L}, we have for $%
\left\vert \beta \right\vert >1$%
\begin{equation*}
\left\vert \beta \right\vert \left\Vert p_{0}\right\Vert _{\mathcal{O}%
}^{2}\leq C\left\vert \beta \right\vert \left\Vert p_{0}\right\Vert _{%
\mathcal{O}}\left( \sqrt{\left\vert \left( \Phi ^{\ast },\Phi \right) _{%
\mathcal{H}}\right\vert }+\left\Vert \Phi ^{\ast }\right\Vert _{\mathcal{H}%
}\right) +\epsilon \left\Vert \Phi \right\Vert _{\mathcal{H}%
}^{2}+C_{\epsilon }\left\Vert \Phi ^{\ast }\right\Vert _{\mathcal{H}}^{2};
\end{equation*}%
and so for $\left\vert \beta \right\vert >1$, we obtain upon division%
\begin{equation*}
\left\Vert p_{0}\right\Vert _{\mathcal{O}}^{2}\leq C\left\Vert
p_{0}\right\Vert _{\mathcal{O}}\left( \sqrt{\left\vert \left( \Phi ^{\ast
},\Phi \right) _{\mathcal{H}}\right\vert }+\left\Vert \Phi ^{\ast
}\right\Vert _{\mathcal{H}}\right) +\epsilon \left\Vert \Phi \right\Vert _{%
\mathcal{H}}^{2}+C_{\epsilon }\left\Vert \Phi ^{\ast }\right\Vert _{\mathcal{%
H}}^{2}.
\end{equation*}%
Applying Young's Inequality with $1/2>\delta >0$, we now obtain in the case $%
\left\vert \beta \right\vert >1$,%
\begin{equation}
(1-\delta )\left\Vert p_{0}\right\Vert _{\mathcal{O}}^{2}\leq \epsilon
\left\Vert \Phi \right\Vert _{\mathcal{H}}^{2}+C_{\delta }\left\vert \left(
\Phi ^{\ast },\Phi \right) _{\mathcal{H}}\right\vert +C_{\epsilon ,\delta
}\left\Vert \Phi ^{\ast }\right\Vert _{\mathcal{H}}^{2}.  \label{5.36}
\end{equation}

\medskip

\textbf{Case II (}$\left\vert \beta \right\vert \leq 1$\textbf{):}\newline
Combining the decomposition (\ref{pre}) with Lemmas \ref{med} and \ref{US-L}%
, we have for the case $\left\vert \beta \right\vert \leq 1$, 
\begin{eqnarray}
\left\Vert p_{0}\right\Vert _{\mathcal{O}}^{2} &=&\left\Vert
q_{0}\right\Vert _{\mathcal{O}}^{2}+\left\Vert c_{0}\right\Vert _{\mathcal{O}%
}^{2}  \label{5.37} \\
&\leq &C\left( \left\Vert \Phi \right\Vert _{\mathcal{H}}\left\Vert \Phi
^{\ast }\right\Vert _{\mathcal{H}}+\left\Vert \Phi ^{\ast }\right\Vert _{%
\mathcal{H}}^{2}+\left\vert \beta \right\vert ^{2}\left\Vert \mathbf{u}%
_{0}\right\Vert _{\mathcal{O}}^{2}\right)  \notag \\
&\leq &C\left( \left\Vert \Phi \right\Vert _{\mathcal{H}}\left\Vert \Phi
^{\ast }\right\Vert _{\mathcal{H}}+\left\Vert \Phi ^{\ast }\right\Vert _{%
\mathcal{H}}^{2}\right) .  \notag
\end{eqnarray}

\medskip

With estimates (\ref{5.36}) and (\ref{5.37}) in hand, we have now upon
invoking $\left\vert ab\right\vert \leq \epsilon a^{2}+C_{\epsilon }b^{2}$
and a rescaling of $\epsilon >0$, that for all $\beta \in \mathbb{R}$,%
\begin{equation}
\left\Vert p_{0}\right\Vert _{\mathcal{O}}^{2}\leq \epsilon \left\Vert \Phi
\right\Vert _{\mathcal{H}}^{2}+C_{\epsilon }\left\Vert \Phi ^{\ast
}\right\Vert _{\mathcal{H}}^{2}.  \label{5.38}
\end{equation}

\bigskip

\emph{The Conclusion of the Proof of Theorem \ref{uniform}:} Combining the
estimates in Lemma \ref{US-L} and (\ref{5.38}), which we have obtained for
the solution $\left[ p_{0},\mathbf{u}_{0},w_{1},w_{2}\right] =\Phi $ of the
resolvent equation (\ref{5.2}) with given data \newline
$\Phi ^{\ast }=\left[ p_{0}^{\ast },\mathbf{u}_{0}^{\ast },w_{1}^{\ast
},w_{2}^{\ast }\right] \in H_{N}^{\bot }$, we have 
\begin{equation*}
\left\Vert \Phi \right\Vert _{\mathcal{H}}^{2}\leq \epsilon \left\Vert \Phi
\right\Vert _{\mathcal{H}}^{2}+C\left[ \left\vert \left( \Phi ^{\ast },\Phi
\right) _{\mathcal{H}}\right\vert +\left\Vert \Phi ^{\ast }\right\Vert _{%
\mathcal{H}}^{2}\right] .
\end{equation*}%
Applying $\left\vert ab\right\vert \leq \epsilon a^{2}+C_{\epsilon }b^{2}$
one last time, and rescaling $\epsilon >0$ sufficiently small, we have
finally the estimate, for given $\beta \in \mathbb{R}$,%
\begin{equation}
\left\Vert (i\beta -\mathcal{A})\Phi ^{\ast }\right\Vert _{\mathcal{H}%
}^{2}\leq C\left\Vert \Phi ^{\ast }\right\Vert _{\mathcal{H}}^{2}.
\label{5.39}
\end{equation}%
Since $\Phi ^{\ast }\in H_{N}^{\bot }$ is arbitrary, (\ref{5.39})
establishes the desired resolvent estimate (\ref{5.1}). The proof of Theorem %
\ref{uniform} is accomplished.

\section{Appendix}

Throughout this paper we have appealed to various results from applied
functional analysis and PDE theory. For the reader's convenience we collect
some of these below.

The following classic lemma provides a useful criterion for the existence of
bounded inverse of a linear operator, and is crucial in the investigation of
the resolvent set of the compressible flow-structure generator:

\begin{lemma}
\label{inverse} (See e.g., Lemma 3.8.18(b) of \cite{HPC}.) Let $L$ be a
linear, closed, operator from Hilbert space $H\rightarrow H$. Then $%
L^{-1}\in \mathcal{L}(H)\Leftrightarrow Range(L)$ is dense in $H$ and there
is a constant $m>0$ such that \newline
$\Vert Lf\Vert \geq m\Vert f\Vert ,~\forall f\in D(L).$
\end{lemma}

In the course of analyzing the linearized compressible flow-structure PDE
model under present consideration, we utilize the following classical
theorem given for nonhomogeneous Stokes problems :

\begin{theorem}
(See Theorem 2.4 and Remark 2.5 of \cite{temam}; see also \cite{consta} and 
\cite{sohr}.) \label{static} With $\Omega $ being a Lipschitz bounded domain
in $\mathbb{R}^{n}$, let data $\left\{ \mathbf{f},g,\mathbf{\phi }\right\}
\in \mathbf{H}^{-1}(\Omega )\times L^{2}(\Omega )\times \mathbf{H}%
^{1/2}(\partial \Omega )$ be given, with $\mathbf{\phi }$ and $g$
furthermore satisfying the compatibility condition%
\begin{equation*}
\int_{\Omega }g(x)dx=\int_{\partial \Omega }\mathbf{\phi }\cdot \nu d\Gamma .
\end{equation*}%
Then there exists $\left\{ \mathbf{u,}p\right\} \in H^{1}(\Omega )\times
L^{2}(\Omega )$ which are solutions of the nonhomogeneous Stokes problem,

\begin{equation*}
\left\{ 
\begin{array}{l}
-\nu \Delta \mathbf{u}+\nabla p=\mathbf{f}~\text{ in }~\Omega \\ 
div~\mathbf{u}=g~\text{ in }~\Omega \\ 
\mathbf{u}=\mathbf{\phi }~\text{ on }~\partial \Omega%
\end{array}%
\right.
\end{equation*}%
$\mathbf{u}$ is unique and $p$ is unique up to the addition of a constant.
Moreover, one has the following estimate:%
\begin{equation*}
\left\Vert \mathbf{u}\right\Vert _{\mathbf{H}^{1}(\Omega )}+\left\Vert
p\right\Vert _{\frac{L^{2}(\Omega )}{\mathbf{R}}}\leq C\left( \left\Vert [%
\mathbf{f},g,\mathbf{\phi }]\right\Vert _{\mathbf{H}^{-1}(\Omega )\times
L^{2}(\Omega )\times \mathbf{H}^{1/2}(\partial \Omega )}\right) .
\end{equation*}
\end{theorem}

We recall next the strong stability resolvent criterion for strongly
continuous semigroups:

\begin{theorem}
\label{stab-st}(see \cite{arendt}.) Let ${\{T(t)\}}_{t\geq 0}$ be a bounded $%
C_{0}$-semigroup with generator $A$ on a reflexive space $X$. Assume that no
eigenvalue of $A$ lies on the imaginary axis. Then if $\sigma (A)\cap i%
\mathbb{R}$ is countable, ${\{T(t)\}}_{t\geq 0}$ is stable.
\end{theorem}

\bigskip

The uniform stability criterion which is used in the proof of our main
result Theorem \ref{uniform} is as follows:

\begin{theorem}
\label{stab-un} (see \cite{huang}, \cite{pruss}.) Let $\left\{
e^{At}\right\} _{t\geq 0}$ be a $C_{0}$-semigroup generated by $A$ in a
Hilbert space which satisfies%
\begin{equation*}
\left\Vert e^{At}\right\Vert \leq K_{0}\text{, for every }t\geq 0\text{, for
some constant }K_{0}>0\text{.}
\end{equation*}%
Then $e^{At}$ decays exponentially iff:

\begin{description}
\item[(i)] $\left\{ \lambda \in \mathbb{C}:\lambda =i\omega \text{, }\omega
\in \mathbb{R}\right\} \subset \rho (A)$;

\item[(ii)] 
\begin{equation*}
\sup_{\omega \in \mathbb{R}}\left\Vert (i\omega -A)^{-1}\right\Vert <\infty .
\end{equation*}
\end{description}
\end{theorem}


\bigskip

\begin{thebibliography}{99}
\bibitem{arendt} W. Arendt and C.J.K Batty, \textquotedblleft Tauberian
theorems and stability of one-parameter semigroups\textquotedblright , \emph{%
Transactions of the American Mathematical Society}, Vol. 306,Number 2, pp.
837-852 (April 1988).

\bibitem{decay} Kagei, Y., \textquotedblleft Decay estimates on solutions of
the linearized compressible Navier-Stokes equation around a parallel flow in
a cylindrical domain, \emph{Kyusha J. Math. }69 (2015), pp. 293-343.

\bibitem{spectral} Aoyama, R. and Kagei, Y., 2016. Spectral properties of
the semigroup for the linearized compressible Navier-Stokes equation around
a parallel in a cylindrical domain. \emph{Advances in Differential Equations}%
, 21(3/4), pp.265--300.

\bibitem{aubin} Aubin, J.P., 2011. \emph{Applied functional analysis} (Vol.
47). John Wiley \& Sons.

\bibitem{george1} Avalos, G. and Bucci, F., 2014. Exponential decay
properties of a mathematical model for a certain flow-structure interaction.
In \emph{New Prospects in Direct, Inverse and Control Problems for Evolution
Equations} (pp. 49--78). Springer International Publishing.

\bibitem{george2} Avalos, G. and Bucci, F., 2015. Rational rates of uniform
decay for strong solutions to a flow-structure PDE system. \emph{Journal of
Differential Equations}, 258(12), pp.4398--4423.

\bibitem{agw} G. Avalos, P. G. Geredeli and J.T. Webster \textquotedblleft
Semigroup Well-posedness of A Linearized, Compressible flow with An Elastic
Boundary\textquotedblright , arxiv.org/abs/1703.10855, and to appear in 
\emph{Discrete and Continuous Dynamical Systems.}

\bibitem{ALT} G. Avalos, R. Triggiani, and I. Lasiecka, Heat-Wave
interaction in 2 or 3 dimensions: optimal decay rates\textquotedblright , 
\emph{Journal of Mathematical Analysis and Applications, }Volume 437, Issue
2, 15 May 2016, Pages 782--815.


\bibitem{bolotin} Bolotin, V.V., 1963. \emph{Nonconservative problems of the
theory of elastic stability}. Macmillan.

\bibitem{buffa2} Buffa, A., Costabel, M. and Sheen, D., 2002. On traces for $%
\mathbf{H}(\text{curl},\Omega )$ in Lipschitz domains. \emph{Journal of
Mathematical Analysis and Applications}, 276(2), pp.845--867.


\bibitem{chorin-marsden} Chorin, A.J. and Marsden, J.E., 1990. \emph{A
mathematical introduction to flow mechanics} (Vol. 3). New York: Springer.

\bibitem{Chu2013-comp} Chueshov, I., 2014. Dynamics of a nonlinear elastic
plate interacting with a linearized compressible viscous flow. \emph{%
Nonlinear Analysis: Theory, Methods \& Applications}, 95, pp.650--665.

\bibitem{Igor-note} Chueshov, I., Personal communication, March 2013.



\bibitem{supersonic} Chueshov, I., Lasiecka, I. and Webster, J.T., 2013.
Evolution semigroups in supersonic -plate interactions. \emph{Journal of
Differential Equations}, 254(4), pp.1741--1773.


\bibitem{consta} P. Constantin and C. Foias, \emph{Navier-Stokes Equations},
The University of Chicago Press, Chicago (1988).

\bibitem{dV} da Veiga, H.B., 1985. \emph{Stationary Motions and
Incompressible Limit for Compressible Viscous flows}, Houston Journal of
Mathematics, Volume 13, No. 4 (1987), pp. 527-544.

\bibitem{dauge} M. Dauge, January 1989. \emph{Stationary Stokes and Navier
Stokes Systems on Two or Three Dimensional Domains with Corners, Part I:
Linearized Equations}, Siam J. Math. Anal., Vol 20, No.1.

\bibitem{deckel} K. Deckelnick, Decay estimates for the compressible
Navier--Stokes equations in unbounded domain, \emph{Math. Z.} 209 (1992),
115--130.

\bibitem{dowell1} E. Dowell, 2004. \emph{A Modern Course in Aeroelasticity}. 
{Kluwer Academic Publishers}.

\bibitem{friedman} Friedman, B. \emph{Principles and Techniques of Applied
Mathematics}, Dover Publications, Inc., New York (1990).

\bibitem{grisvard} Grisvard, P., 2011. \emph{Elliptic problems in nonsmooth
domains}. Society for Industrial and Applied Mathematics.

\bibitem{galdi} G.P. Galdi, \emph{An introduction to the mathematical theory
of the Navier-Stokes equations Volume I}, Springer Tracts in Natural
Philosophy, Springer-Verlag, New York (1991).

\bibitem{huang} Huang, F.L. Characteristic conditions for exponential
stability of linear dynamical systems in Hilbert spaces, \emph{Ann. Differ.
Equ.}, 1(1), pp. 43-53 (1985).

\bibitem{HPC} Hutson, V., Pym, J. and Cloud, M., 2005. \emph{Applications of
functional analysis and operator theory, Second Edition}, Elsevier, New
York, 2005.

\bibitem{j-k} D.S. Jerrison and C.E. Kenig, The Neumann problem on Lipschitz
domains, Bulletin of the AMS, Volume 4, Number 2 (March 1981), pp. 203-207.

\bibitem{kato} Kato, T., 2013. \emph{Perturbation theory for linear operators%
} (Vol. 132). Springer Science \& Business Media.

\bibitem{kagei} Y. Kagei and T. Kobayashi, Asymptotic behavior of solutions
to the compressible Navier--Stokes equations on the half space, \emph{Arch.
Ration. Mech. Anal. }177 (2005), 231--330.

\bibitem{kesavan} Kesavan, S., 1989. \emph{Topics in functional analysis and
applications}.

\bibitem{kob} T. Kobayashi. Some estimates of solutions for the equations of
motion of compressible viscous fluid in an exterior domain in $\mathbb{R}%
^{3} $, \emph{J. Differential Equations} 184 (2002), 587--619.




\bibitem{Mc} McLean, W.C.H., 2000. \emph{Strongly elliptic systems and
boundary integral equations}. Cambridge university press.

\bibitem{necas} Ne\v{c}as, 2012. Direct Methods in the Theory of Elliptic
Equations (translated by Gerard Tronel and Alois Kufner), Springer, New York.

\bibitem{pazy} Pazy, A., 2012. \emph{Semigroups of linear operators and
applications to partial differential equations} (Vol. 44). Springer Science
\& Business Media.

\bibitem{pruss} Pr$\ddot{u}$ss, J., August 1984. \emph{On the spectrum of $%
C_{0}$-Semigroups}, Transactions of the American Mathematical Society,
Volume 284, Number 2.

%

\bibitem{sohr} H. Sohr, \emph{The Navier-Stokes Equations, An Elementary
Functional Analytical Approach}, Birkh\"{a}user Verlag, Boston (2001).

\bibitem{temam} Temam, R., 1977. \emph{Navier-Stokes Equations Theory and
Numerical Analysis}, North-Holland Publishing Company,
Amsterdam-NewYork-Oxford.

\bibitem{valli} Valli, A., 1987. On the existence of stationary solutions to
compressible Navier-Stokes equations. In \emph{Annales de l'IHP Analyse non
lin\'{e}aire} (Vol. 4, No. 1, pp. 99--113).

\end{thebibliography}
\end{document}